\documentclass{ceb-l}

\usepackage{amsmath}
\usepackage{amssymb}
\usepackage{euscript}

\usepackage[pdftex]{graphicx}
\newcommand{\pix}[3]{\begin{figure}[tbp]
\centering{\includegraphics[width=#1\textwidth]{#2}}
\caption{#3}
\end{figure}}

\newtheorem{theorem}{Theorem}[section]
\newtheorem{lemma}[theorem]{Lemma}
\newtheorem{prop}[theorem]{Proposition}
\newtheorem{coro}[theorem]{Corollary}

\theoremstyle{definition}

\newtheorem{example}[theorem]{Example}
\newtheorem{exercise}[theorem]{Exercise}

\theoremstyle{remark}

\numberwithin{equation}{section}

\newcommand{\showon}{\begin{eqnarray*}}
\newcommand{\showoff}{\end{eqnarray*}}
\newcommand{\drop}{\smallsetminus}
\newcommand{\goesto}{\rightarrow}
\newcommand{\diag}{\mathrm{diag}}
\newcommand{\zero}{\boldsymbol{0}}
\newcommand{\one}{\boldsymbol{1}}
\newcommand{\none}{\varnothing}

\newcommand{\CAP}{\mathrm{cap}}
\newcommand{\tr}{\mathrm{tr}}
\newcommand{\Alg}{\mathfrak{A}}
\renewcommand{\AA}{\mathbf{A}}
\newcommand{\A}{\EuScript{A}}

\newcommand{\D}{\EuScript{D}}
\newcommand{\E}{\EuScript{E}}
\newcommand{\F}{\EuScript{F}}
\newcommand{\J}{\EuScript{J}}
\newcommand{\EE}{\mathfrak{E}}
\newcommand{\Pol}{\mathrm{Pol}}
\renewcommand{\Pr}{\mathrm{Pr}}
\newcommand{\NN}{\mathbb{N}}
\newcommand{\ZZ}{\mathbb{Z}}
\newcommand{\RR}{\mathbb{R}}
\newcommand{\CC}{\mathbb{C}}
\newcommand{\KK}{\mathbb{K}}
\newcommand{\MM}{\EuScript{M}}
\newcommand{\G}{\mathrm{G}}
\newcommand{\HH}{\EuScript{H}}
\newcommand{\Hbar}{\overline{\HH}}
\renewcommand{\L}{\EuScript{L}}
\renewcommand{\i}{\mathrm{i}}
\newcommand{\e}{\mathrm{e}}
\renewcommand{\Im}{\mathrm{Im}}
\renewcommand{\Re}{\mathrm{Re}}
\renewcommand{\S}{\EuScript{S}}
\newcommand{\Stab}{\mathfrak{S}}
\newcommand{\StabR}{\Stab_\RR}
\newcommand{\ma}{\mathsf{MA}}
\newcommand{\T}{\mathsf{T}}
\newcommand{\Wr}{\mathrm{W}}
\renewcommand{\u}{\mathbf{u}}
\renewcommand{\v}{\mathbf{v}}
\newcommand{\w}{\mathbf{w}}
\newcommand{\x}{\mathbf{x}}
\newcommand{\y}{\mathbf{y}}
\newcommand{\z}{\mathbf{z}}
\renewcommand{\a}{\mathbf{a}}
\renewcommand{\b}{\mathbf{b}}
\renewcommand{\c}{\mathbf{c}}

\newcommand{\fhat}{\widehat{f}}
\newcommand{\ghat}{\widehat{g}}
\renewcommand{\aa}{\mathfrak{a}}
\newcommand{\ff}{\mathfrak{f}}
\newcommand{\pp}{\mathfrak{t}}
\newcommand{\del}{\boldsymbol{\partial}}
\newcommand{\Disc}{\mathrm{Disc}}
\newcommand{\Det}{\mathrm{Det}}
\newcommand{\per}{\mathrm{per}}
\newcommand{\supp}{\mathrm{supp}}

\begin{document}

\title[Multivariate stable polynomials]{Multivariate stable polynomials:\\
theory and applications}


\author{David G. Wagner}
\address{Department of Combinatorics and Optimization,
University of Waterloo, Waterloo, Ontario, Canada N2L 3G1}
\email{dgwagner@math.uwaterloo.ca}
\thanks{Research supported by NSERC Discovery Grant OGP0105392.}


\dedicatory{In memoriam Julius Borcea.}

\begin{abstract}
Univariate polynomials with only real roots -- while special -- do occur often
enough that their properties can lead to interesting conclusions in diverse areas.
Due mainly to the recent work of two young mathematicians, Julius Borcea and
Petter Br\"and\'en, a very successful multivariate generalization of this method has
been developed.  The first part of this paper surveys some of the main
results of this theory of ``multivariate stable'' polynomials -- the most central
of these results is the characterization of linear transformations preserving 
stability of polynomials.  The second part presents various applications of this 
theory in complex analysis, matrix theory, probability and statistical mechanics,
and combinatorics.
\end{abstract}

\maketitle

\section{Introduction.}

I have been asked by the AMS to survey the recent work of Julius Borcea
and Petter Br\"and\'en on their multivariate generalization of the theory
of univariate polynomials with only real roots, and its applications.
It is exciting work -- elementary but subtle, and with spectacular consequences.
Borcea and Br\"and\'en take center stage but there are many other actors,
many of whom I am unable to mention in this brief treatment.  Notably, Leonid
Gurvits provides a transparent proof of a vast generalization of the famous van
der Waerden Conjecture.

Space is limited and I have been advised to use ``Bourbaki style'',
and so this is an account of the essentials of the theory and a 
few of its applications, with complete proofs as far as possible.  Some 
relatively straightforward arguments have been left as exercises to 
engage the reader, and some more specialized topics are merely sketched 
or even omitted.  For the full story and the history and context of the subject
one must go to the references cited, the references they cite, and so on.
The introduction of \cite{BB4}, in particular, gives a good account 
of the genesis of the theory.

Here is a brief summary of the contents.  Section 2 introduces stable
polynomials, gives some examples, presents their elementary properties,
and develops multivariate generalizations of two classical univariate
results:\ the Hermite-Kakeya-Obreschkoff and Hermite-Biehler Theorems.
We also state the P\'olya-Schur Theorem characterizing ``multiplier sequences'',
as this provides an inspiration for much of the multivariate theory. 
Section 3 restricts attention to multiaffine stable polynomials:\
we present a characterization of multiaffine real stable polynomials
by means of parameterized quadratic inequalities, and characterize those
linear transformations which take multiaffine stable polynomials
to stable polynomials.  In Section 4 we use parts of the forgoing for
Borcea and Br\"and\'en's splendid proof of the Grace-Walsh-Szeg\H{o}
Coincidence Theorem.  In Section 5, the Grace-Walsh-Szeg\H{o} Theorem is
used to extend the results of Section 3 from multiaffine to arbitrary stable
polynomials.  This culminates in an amazing multivariate generalization
of the P\'olya-Schur Theorem, the proof of which requires the 
development of a multivariate extension of the Szasz Principle (which
is omitted, regretfully,  for lack of space). 
Section 6 presents Borcea and Br\"and\'en's resolution of some matrix-theoretic
conjectures of Johnson.  Section 7 presents the derivation by Borcea, Br\"and\'en,
and Liggett of negative association inequalities for the symmetric exclusion process,
a fundamental model in probability and statistical mechanics.
Section 8 presents Gurvits's sweeping generalization of the van
der Waerden Conjecture.  Finally, Section 9 briefly mentions a few
further topics that could not be included fully for lack of space.

I thank Petter Br\"and\'en kindly for his helpful comments on preliminary 
drafts of this paper.

\section{Stable polynomials.}

We use the following shorthand notation for multivariate polynomials.
Let $[m]=\{1,2,...,m\}$, let $\x=(x_1,...,x_m)$ be a sequence of
indeterminates, and let $\CC[\x]$ be the ring of complex polynomials
in the indeterminates $\x$.  For a function $\alpha:[m]\goesto\NN$, let
$\x^\alpha = x_1^{\alpha(1)}\cdots x_m^{\alpha(m)}$ be the corresponding
monomial.  For $S\subseteq[m]$ we also let $\x^S = \prod_{i\in S} x_i$.
Similarly, for $i\in[m]$ let $\partial_i = \partial/\partial x_i$,
let $\del = (\partial_1,...,\partial_m)$, let
$\del^\alpha = \partial_1^{\alpha(1)}\cdots \partial_m^{\alpha(m)}$
and let $\del^S = \prod_{i\in S} \partial_i$.
The constant functions on $[m]$ with images $0$ or $1$ are denoted by
$\zero$ and $\one$, respectively.
The $\x$ indeterminates are always indexed by $[m]$.

Let $\HH=\{z\in\CC:\ \Im(z)>0\}$ denote the open upper half of the complex plane,
and $\Hbar$ the closure of $\HH$ in $\CC$.  A polynomial $f\in\CC[\x]$ is
\emph{stable} provided that either $f\equiv 0$ identically, or
whenever $\z=(z_1,...,z_m)\in\HH^m$ then $f(\z)\neq 0$.  We use $\Stab[\x]$
to denote the set of stable polynomials in $\CC[\x]$, and
$\StabR[\x]=\Stab[\x]\cap\RR[\x]$ for the set of \emph{real stable} polynomials
in $\RR[\x]$. (Borcea and Br\"and\'en do not consider the zero polynomial
to be stable, but I find the above convention more convenient.)

We rely on the following essential fact at several points.\\

\noindent
\textbf{Hurwitz's Theorem}\ (Theorem 1.3.8 of \cite{RS})\textbf{.}\
\emph{Let $\Omega\subseteq\CC^m$ be a connected
open set, and let $(f_n:\ n\in\NN)$ be a sequence of functions, each
analytic and nonvanishing on $\Omega$, which converges to a limit $f$ uniformly
on compact subsets of $\Omega$.  Then $f$ is either nonvanishing on $\Omega$
or identically zero.}\\

Consequently, a polynomial obtained as the limit of a convergent sequence of
stable polynomials is itself stable.

\subsection{Examples.}

\begin{prop}[Proposition 2.4 of \cite{BB1}]
For $i\in[m]$, let $A_i$ be an $n$-by-$n$ matrix and let $x_i$ 
be an indeterminate, and let $B$ be an $n$-by-$n$ matrix. If $A_i$ is
positive semidefinite for all $i\in[m]$ and $B$ is Hermitian then
$$
f(\x) = \det(x_1 A_1 + x_2 A_2 + \cdots + x_m A_m + B)
$$
is real stable.
\end{prop}
\begin{proof}
Let $\overline{f}$ denote the coefficientwise complex conjugate of $f$.
Since $\overline{A_i}=A_i^\T$ for all $i\in[m]$, and $\overline{B}=B^\T$,
it follows that $\overline{f}=f$, so that $f\in\RR[\x]$.  By Hurwitz's Theorem
and a routine perturbation argument,
it suffices to prove that $f$ is stable when each $A_i$ is positive definite.  
Consider any $\z=\a+\i\b\in\HH^m$, with $\a,\b\in\RR^m$ and $b_i>0$ for all
$i\in[m]$ (abbreviated to $\b>\zero$).  Now $Q=\sum_{i=1}^m b_i A_i$ is
positive definite, and hence has a
positive definite square-root $Q^{1/2}$.  Also note that
$H=\sum_{i=1}^m a_i A_i + B$ is Hermitian, and that
$$
f(\z) = \det(Q)\det(\i I + Q^{-1/2}HQ^{-1/2}).
$$
Since $\det(Q)\neq 0$, if $f(\z)=0$ then $-\i$ is an eigenvalue of
$Q^{-1/2}HQ^{-1/2}$, contradicting the fact that this matrix is Hermitian.
Thus, $f(\z)\neq 0$ for all $\z\in\HH^m$.  That is, $f$ is stable.
\end{proof}

\begin{coro}
Let $Q$ be an $n$-by-$m$ complex matrix, and let
$X=\diag(x_1,...,x_m)$ be a diagonal matrix of indeterminates.
Then $f(\x)=\det(QXQ^\dagger)$ is real stable.
\end{coro}
\begin{proof}
Let $Q=(q_{ij})$, and for $j\in[m]$ let $A_j$ denote
the $n$-by-$n$ matrix with $hi$-th entry $q_{hj}\overline{q}_{ij}$.  That is,
$A_j=Q_j Q_j^\dagger$ in which $Q_j$ denotes the $j$-th column of $Q$.
Since each $A_j$ is positive semidefinite and
$QXQ^\dagger = x_1 A_1 + \cdots + x_m A_m$,
the conclusion follows directly from Proposition 2.1.
\end{proof}

\subsection{Elementary properties.}

The following simple observation often allows multivariate problems to
be reduced to univariate ones, as will be seen. 
\begin{lemma}
A polynomial $f\in\CC[\x]$ is stable if and only if for all
$\a,\b\in\RR^m$ with $\b>\zero$, $f(\a + \b t)$ is stable in $\Stab[t]$.
\end{lemma}
\begin{proof}
Since $\HH^{m} = \{\a+\b t:\ \a,\b\in\RR^{m},\ \b>\zero,\
\mathrm{and}\ t\in\HH\}$, the result follows.
\end{proof}

For $f\in\CC[\x]$ and $i\in[m]$, let $\deg_i(f)$ denote the degree of $x_i$
in $f$.
\begin{lemma}
These operations preserve stability of polynomials in $\CC[\x]$.\\
\textup{(a)}\ \textbf{Permutation:}\ for any permutation $\sigma:[m]\goesto[m]$,
$f \mapsto f(x_{\sigma(1)},...,x_{\sigma(m)})$.\\
\textup{(b)}\ \textbf{Scaling:}\ for $c\in\CC$ and $\a\in\RR^{m}$ with
$\a>\zero$, $f \mapsto cf(a_{1}x_{1},\ldots,a_{m}x_{m})$.\\
\textup{(c)}\ \textbf{Diagonalization:}\ for $\{i,j\}\subseteq[m]$,
$f \mapsto f(\x)|_{x_i=x_j}$.\\
\textup{(d)}\ \textbf{Specialization:}\ for $a\in\Hbar$,
$f \mapsto f(a, x_{2},\ldots,x_{m})$.\\
\textup{(e)}\ \textbf{Inversion:}\ if $\deg_1(f)=d$,
$f \mapsto x_{1}^{d}f(-x_{1}^{-1},x_{2},\ldots,x_{m})$.\\
\textup{(f)}\ \textbf{Differentiation (or ``Contraction''):}\
$f \mapsto \partial_1 f(\x)$.
\end{lemma}
\begin{proof}
Parts (a,b,c) are clear.  Part (d) is also clear in the case that
$\Im(a)>0$.  For $a\in\RR$ apply part (d) with values in the sequence
$(a + \i 2^{-n}:\ n\in\NN)$, and then apply Hurwitz's Theorem to the
limit as $n\goesto\infty$.  Part (e) follows from the fact that
$\HH$ is invariant under the operation $z\mapsto -z^{-1}$.  For part (f),
let $d = \deg_{1}(f)$, and consider the sequence $f_{n} = n^{-d} 
f(nx_{1},x_{2},\ldots,x_{m})$ for all $n\geq 1$.  Each $f_{n}$
is stable and the sequence converges to a polynomial, so the limit is stable.
Since $\deg_1(f)=d$, this limit is not identically zero.  This implies that
for all $z_2,...,z_m\in\HH$, the polynomial $g(x)=f(x,z_2,...,z_m)\in\CC[x]$
has degree $d$.  Clearly $g'(x)=\partial_1 f(x,z_2,...,z_m)$.  Let
$\xi_1,...,\xi_d$ be the roots of $g(x)$, so that
$g(x)=c\prod_{h=1}^d (x-\xi_h)$ for some $c\in\CC$.
Since $f$ is stable, $\Im(\xi_h)\leq 0$ for all $h\in[d]$.  Now
$$
\frac{g'(x)}{g(x)}=\frac{d}{dx}\log g(x)=\sum_{h=1}^d\frac{1}{x-\xi_h}.
$$
If $\Im(z)>0$ then $\Im(1/(z-\xi_h))<0$ for all $h\in[d]$, so that
$g'(z)\neq 0$.  Thus, if $\z\in\HH^m$ then $\partial_1 f(\z)\neq 0$.
That is, $\partial_1 f$ is stable.
\end{proof}
\noindent
Of course, by permutation, parts (d,e,f) of Lemma 2.4 apply for any index
$i\in[m]$ as well (not just $i=1$).  Part (f) is essentially the Gauss-Lucas
Theorem:\ the roots of $g'(x)$ lie in the convex hull of the roots of $g(x)$.

\subsection{Univariate stable polynomials.}

A nonzero univariate polynomial is real stable if and only if
it has only real roots.  Let $f$ and $g$ be two such polynomials,
let $\xi_{1}\leq\xi_{2}\leq\cdots\leq\xi_{k}$ be the roots of $f$,
and let $\theta_{1}\leq\theta_{2}\leq\cdots\leq\theta_{\ell}$ be the 
roots of $g$.  These roots are \emph{interlaced} if they are ordered so that 
$\xi_{1}\leq\theta_{1}\leq\xi_{2}\leq\theta_{2}\leq\cdots$ or
$\theta_{1}\leq\xi_{1}\leq\theta_{2}\leq\xi_{2}\leq\cdots$.
For each $i\in[\ell]$, let $\ghat_{i} = g/(x-\theta_{i})$.  If $\deg f 
\leq \deg g$ and the roots of $g$ are simple,
then there is a unique $(a,b_{1},\ldots,b_{\ell})\in\RR^{\ell+1}$
such that
$$
f = ag + b_{1} \ghat_{1} + \cdots + b_{\ell} \ghat_{\ell}.
$$
\begin{exercise}
Let $f,g\in\StabR[x]$ be nonzero and such that $fg$ has only simple roots,
let $\deg f \leq \deg g$,  and let $\theta_{1}<\cdots<\theta_{\ell}$ be
the roots of $g$.  The following are equivalent:\\
\textup{(a)}\  The roots of $f$ and $g$ are interlaced.\\
\textup{(b)}\ The sequence  $f(\theta_{1}),f(\theta_{2}),\ldots,f(\theta_{\ell})$
alternates in sign (strictly).\\ 
\textup{(c)}\ In $f = ag + \sum_{i=1}^{\ell} b_{i} \ghat_{i}$, all of
$b_{1},\ldots,b_{\ell}$ have the same sign (and are nonzero).
\end{exercise}

The\emph{Wronskian} of $f,g\in \CC[x]$ is $\Wr[f,g]=f'\cdot g-f\cdot g'$.
If $f = ag + \sum_{i=1}^{\ell} b_{i} \ghat_{i}$ as in Exercise 2.5 then
$$
\frac{\Wr[f,g]}{g^{2}} = \frac{d}{dx}\left(\frac{f}{g}\right) =
\sum_{i=1}^{\ell} \frac{-b_{i}}{(x-\theta_{i})^{2}}.
$$
It follows that if $f$ and $g$ are as in Exercise 2.5(a) then
$\Wr[f,g]$ is either positive for all real $x$, or negative for
all real $x$.   Since $\Wr[g,f]=-\Wr[f,g]$ the condition that
$\deg f\leq \deg g$ is immaterial.  Any pair $f,g$ with interlacing 
roots can be approximated arbitrarily closely by such a pair with all
roots of $fg$ simple.  It follows that for any pair $f,g$ with interlacing 
roots, the Wronskian $\Wr[f,g]$ is either nonnegative on all of $\RR$ 
or nonpositive on all of $\RR$.

Nonzero univariate polynomials $f,g\in\StabR[x]$ are in \emph{proper position}, 
denoted by $f\ll g$, if $\Wr[f,g]\leq 0$ on all of $\RR$.
For convenience we also let $0\ll f$ and $f\ll 0$ for any 
$f\in\StabR[x]$;\ in particular $0\ll 0$.

\begin{exercise}
Let $f,g\in\StabR[x]$ be real stable.  Then $f \ll g$ and $g\ll f$
if and only if $cf=dg$ for some $c,d\in\RR$ not both zero.
\end{exercise}

\noindent
\textbf{Hermite-Kakeya-Obreschkoff (HKO) Theorem}
(Theorem 6.3.8 of \cite{RS})\textbf{.}\ \emph{
Let $f,g\in\RR[x]$.  Then $af+bg\in\StabR[x]$ for all $a,b\in\RR$
if and only if $f,g\in\StabR[x]$  and either $f\ll g$ or $g\ll f$.
}\\

\noindent
\textbf{Hermite-Biehler (HB) Theorem}
(Theorem 6.3.4 of \cite{RS})\textbf{.}\ \emph{
Let $f,g\in\RR[x]$.  Then $g+\i f\in\Stab[x]$ if and only if 
$f,g\in\StabR[x]$ and $f\ll g$.
}\\

\begin{proof}[Proofs of HKO and HB]
It suffices to prove these when $fg$ has only simple roots.

For HKO we can assume that $\deg(f)\leq\deg(g)$.  Exercise 2.5
shows that if the roots of $f$ and $g$ are interlaced then for all $a,b\in\RR$,
the roots of $g$ and $af+bg$ are interlaced, so that $af+bg$ is real stable.
The converse is trivial if $cf=dg$ for some $c,d\in\RR$ not both zero, so
assume otherwise.  From the hypothesis, both $f$ and $g$ are real stable.
If there are $z_0,z_1\in\HH$ for which $\Im(f(z_0)/g(z_0))<0$ and
$\Im(f(z_1)/g(z_1))>0$, then for some $\lambda\in[0,1]$ the number
$z_\lambda=(1-\lambda)z_0+\lambda z_1$ is such that
$\Im(f(z_\lambda)/g(z_\lambda))=0$.  Thus $f(z_\lambda)-a g(z_\lambda)=0$
for some real number $a\in\RR$.  Since $f-ag$ is stable (by hypothesis)
and $z_\lambda\in\HH$, this implies that $f-ag\equiv 0$, a contradiction.
Thus $\Im(f(z)/g(z))$ does not change sign for $z\in\HH$.  This implies
Exercise 2.5(c):\  all the $b_i$ have the same sign
(consider $f/g$ at the points $\theta_i+\i \epsilon$ for $\epsilon>0$
approaching $0$).  Thus, the roots of $f$ and $g$ are interlaced.

For HB, let $p=g+\i f$.  Considering $\i p=-f+\i g$ if necessary, we can
assume that $\deg f\leq \deg g$.  If $f\ll g$ then Exercise 2.5(c)
implies that $\Im(f(z)/g(z))\leq 0$ for all $z\in\HH$, so that
$g+\i f$ is stable.  For the converse, let $p(x)=c\prod_{i=1}^d(x-\xi_i)$,
so that $\Im(\xi_i)\leq 0$ for all $i\in[d]$.
Now $|z-\xi_i|\geq |\overline{z}-\xi_i|$ for all $z\in\HH$ and $i\in[d]$,
so that $|p(z)|\geq|p(\overline{z})|$ for all $z\in\HH$.
For any $z\in\HH$ with $f(z)\neq 0$ we have
$$
\left|\frac{g(z)}{f(z)}+\i\right| \geq
\left|\frac{g(\overline{z})}{f(\overline{z})}+\i\right| =
\left|\frac{g(z)}{f(z)}-\i\right|,
$$
and it follows that $\Im(g(z)/f(z))\geq 0$ for all $z\in\HH$ with $f(z)\neq 0$.
Since $fg$ has simple roots it follows that $g(x)+y f(x)$ is stable in
$\Stab[x,y]$.  By contraction and specialization, both $f$ and $g$ are real
stable.  By scaling and specialization, $af+bg$ is stable for all $a,b\in\RR$.
By HKO, the roots of $f$ and $g$ are interlaced.  Since $\Im(f(z)/g(z))\leq 0$
for all $z\in\HH$, all the $b_i$ in Exercise 2.5(c) are positive, so
that $\Wr[f,g]$ is negative on all of $\RR$:\ that is $f\ll g$.
\end{proof}

For $\lambda:\NN\goesto\RR$, let $T_{\lambda}:\RR[x]\goesto\RR[x]$ 
be the linear transformation defined by $T_{\lambda}(x^{n}) = 
\lambda(n) x^{n}$ and linear extension.  A \emph{multiplier sequence
(of the first kind)} is such a $\lambda$ for which $T_{\lambda}(f)$
is real stable whenever $f$ is real stable.  P\'olya and Schur 
characterized multiplier sequences as follows.\\

\noindent
\textbf{P\'olya-Schur Theorem} (Theorem 1.7 of \cite{BB4})\textbf{.}\ \emph{
Let $\lambda:\NN\goesto\RR$.  The following are equivalent:\\
\textup{(a)}\ $\lambda$ is a multiplier sequence.\\
\textup{(b)}\ $F_{\lambda}(x)=\sum_{n=0}^{\infty} \lambda(n) x^{n}/n!$ is an entire 
function which is the limit, uniformly on compact sets, of real
stable polynomials with all roots of the same sign.\\
\textup{(c)}\ Either $F_{\lambda}(x)$ or $F_{\lambda}(-x)$ has the form
$$
Cx^{n}\e^{ax}\prod_{j=1}^{\infty}(1+\alpha_{j}x),
$$
in which $C\in\RR$, $n\in\NN$, $a\geq 0$, all $\alpha_{j}\geq 0$,
and $\sum_{j=1}^{\infty}\alpha_{j}$ is finite.\\
\textup{(d)}\ For all $n\in\NN$, the polynomial 
$T_{\lambda}((1+x)^{n})$ is real stable with all roots of the same 
sign.
}\\

One of the main results of Borcea and Br\"and\'en's theory is a
great generalization of the P\'olya-Schur Theorem -- a 
characterization of \underline{all} \emph{stability preservers}: linear
transformations $T:\CC[\x]\goesto\CC[\x]$ such that $T(f)$ is stable
whenever $f$ is stable.  (Also the analogous characterization of real
stability preservers.)  This is discussed in some detail in Section 5.3.

\subsection{Multivariate analogues of the HKO and HB Theorems.}

By analogy with the univariate HB Theorem, polynomials
$f,g\in\RR[\x]$ are said to be in \emph{proper position}, denoted by
$f\ll g$, when $g+\i f\in\Stab[\x]$.  (As will be seen, this implies
that $f,g\in\StabR[\x]$.)  Thus, the multivariate analogue of
Hermite-Biehler is a \underline{definition}, not a theorem.

\begin{prop}[Lemma 1.8 and Remark 1.3 of \cite{BB5}]
Let $f,g\in\CC[\x]$.\\
\textup{(a)}\ If $f,g\in\RR[\x]$ then $f \ll g$
if and only if $g+yf\in\StabR[\x,y]$.\\
\textup{(b)}\ If $0\not\equiv f\in\Stab[\x]$ then
$g+yf\in\Stab[\x,y]$ if and only if for all $\z\in\HH^m$,
$$
\Im\left(\frac{g(\z)}{f(\z)}\right)\geq 0.
$$
\end{prop}
\begin{proof}
If $g+yf\in\StabR[\x,y]$ then $g+\i f\in\Stab[\x]$, by specialization.  
Conversely, assume that $h=g+\i f\in\Stab[\x]$ with $f,g\in\RR[\x]$,
and let $z=a+\i b$ with $a,b\in\RR$ and $b>0$.  By Lemma 2.3, for all
$\a,\b\in\RR^{m}$ with $\b>\zero$ we have $h(\a+\b t)\in\Stab[t]$.
By HB, $\fhat(t)=f(\a+\b t)$ and $\ghat(t)=g(\a+\b t)$ are
such that $\fhat\ll \ghat$.  By HKO, $c\fhat+d\ghat\in\StabR[t]$
for all $c,d\in\RR$.  By HKO again, the roots of
$b\fhat$ and of $\ghat+a\fhat$ are interlaced.  Since
$W[b\fhat,\ghat+a\fhat] = b W[\fhat,\ghat]\leq 0$ on $\RR$,
it follows that $b\fhat\ll \ghat+a\fhat$.  Finally, by HB again,
$\ghat + (a+\i b)\fhat\in\Stab[t]$.  Since this holds for all
$\a,\b\in\RR^{m}$ with $\b>\zero$, Lemma 2.3 implies that
$g+(a+\i b)f\in\Stab[\x]$.  Since this holds for all
$a,b\in\RR$ with $b>0$, $g+yf\in\Stab[\x,y]$.  This proves part (a).

For part (b), first let $g+yf$ be stable. By specialization, $g$ is also stable.
If $g\equiv 0$ then there is nothing
to prove.  Otherwise, consider any $\z\in\HH^{m}$, so that $f(\z)\neq 0$ and
$g(\z)\neq 0$.  There is a unique solution $z\in\CC$ to $g(\z)+zf(\z)=0$,
and since $g+yf$ is stable, $\Im(z)\leq 0$.  Hence,
$\Im(g(\z)/f(\z))=\Im(-z)\geq 0$.  This argument can be reversed to prove
the converse implication.
\end{proof}

\begin{exercise}[Corollary 2.4 of \cite{BB4}]
$\Stab[\x]=\{g + \i f:\ f,g\in \StabR[\x]\ \mathrm{and}\ f\ll g\}$.
\end{exercise}

Here is the multivariate HKO Theorem of Borcea and Br\"and\'en.

\begin{theorem}[Theorem 1.6 of \cite{BB4}]
Let $f,g\in\RR[\x]$.  Then $af+bg\in\StabR[\x]$ for all $a,b\in\RR$
if and only if $f,g\in\StabR[\x]$ and either $f\ll g$ or $g\ll f$.
\end{theorem}
\begin{proof}
First assume that $f\ll g$, and let $a,b\in\RR$ with $b>0$.
By Proposition 2.7(a), $g+yf\in\StabR[\x,y]$.  By scaling and specialization,
$bg+(a+\i)f \in\Stab[\x]$.  By Proposition 2.7(a) again, $f\ll (af+bg)$.  Thus
$af+bg\in\StabR[\x]$ for all $a,b\in\RR$.  The case that $g\ll f$ is 
similar.

Conversely, assume that $af+bg\in\StabR[\x]$ for all $a,b\in\RR$.
Let $\a,\b\in\RR^{m}$ with $\b>\zero$, and let $\fhat(t)=f(\a+\b t)$ and
$\ghat(t)=g(\a+\b t)$.  By Lemma 2.3, $a\fhat +b\ghat\in\StabR[t]$
for all $a,b\in\RR$.  By HKO, for each
$\a,\b\in\RR^{m}$ with $\b>\zero$, either $\fhat\ll \ghat$
or $\ghat\ll\fhat$.  

If $\fhat\ll \ghat$ for all $\a,\b\in\RR^{m}$ with $\b>\zero$, then
by HB, $\ghat + \i \fhat \in\Stab[t]$ for all
$\a,\b\in\RR^{m}$ with $\b>\zero$.  Thus $g + \i f\in\Stab[\x]$ by
Lemma 2.3, which is to say that $f\ll g$ (by definition).  Similarly, if
$\ghat\ll \fhat$ for all $\a,\b\in\RR^{m}$ with $\b>\zero$ then $g\ll 
f$.

It remains to consider the case that
$f(\a_{0}+\b_{0} t) \ll g(\a_{0}+\b_{0} t)$ for 
some $\a_{0},\b_{0}\in\RR^{m}$ with $\b_{0}>\zero$, and 
$g(\a_{1}+\b_{1} t) \ll f(\a_{1}+\b_{1} t)$ for 
another $\a_{1},\b_{1}\in\RR^{m}$ with $\b_{1}>\zero$.
For $0\leq\lambda\leq 1$, let
$\a_{\lambda} = (1-\lambda)\a_{0}+\lambda\a_{1}$ and 
$\b_{\lambda} = (1-\lambda)\b_{0}+\lambda\b_{1}$.
Since roots of polynomials move continuously as the coefficients 
are varied continuously, there is a value $0\leq \lambda\leq 1$ for 
which both $f(\a_{\lambda}+\b_{\lambda} t) \ll g(\a_{\lambda}+\b_{\lambda} t)$
and $g(\a_{\lambda}+\b_{\lambda} t) \ll f(\a_{\lambda}+\b_{\lambda} t)$.
From Exercise 2.6, it follows that
$cf(\a_{\lambda}+\b_{\lambda} t) = dg(\a_{\lambda}+\b_{\lambda} t)$
for some $c,d\in\RR$ not both zero.  Now $h=cf-dg\in\Stab[\x]$ by hypothesis,
and since $h(\a_{\lambda}+\b_{\lambda} t)\equiv 0$ identically, it follows 
that $h(\a_{\lambda}+\i\b_{\lambda})=0$.  Since $\b_{\lambda}>\zero$ 
and $h$ is stable, this implies that $h\equiv 0$, so that $cf=dg$ in $\Stab[\x]$.
In this case, both $f\ll g$ and $g\ll f$ hold.
\end{proof}

For $f,g\in\CC[\x]$ and $i\in[m]$, let
$\Wr_i[f,g] = \partial_i f \cdot g - f \cdot \partial_i g$
be the \emph{$i$-th Wronskian} of the pair $(f,g)$.
\begin{coro}[Theorem 1.9 of \cite{BB5}]
Let $f,g\in\RR[\x]$.  The following are equivalent:\\
\textup{(a)}\ $g+\i f$ is stable in $\Stab[\x]$, that is $f\ll g$;\\
\textup{(b)}\ $g+ y f$ is real stable in $\StabR[\x,y]$;\\
\textup{(c)}\ $af+bg\in\StabR[\x]$ for all $a,b\in\RR$, and
$\Wr_i[f,g](\a)\leq 0$ for all $i\in[m]$ and $\a\in\RR^m$.
\end{coro}
\begin{proof}
Proposition 2.7(a) shows that (a) and (b) are equivalent.

If (a) holds then Theorem 2.9 implies that $af+bg\in\StabR[\x]$ for all
$a,b\in\RR$.  To prove the rest of (c), let $i\in[m]$ and $\a\in\RR^m$,
and let $\delta_i\in\RR^m$ be the unit vector with a one in the $i$-th position.
Since $f\ll g$, for any $\b\in\RR^m$ with $\b>\zero$ we have
$f(\a+(\b+\delta_i)t) \ll g(\a+(\b+\delta_i)t)$ in $\StabR[t]$, from
Proposition 2.7(a) and Lemma 2.3.  By the Wronskian condition for univariate
polynomials in proper position,
$$
\Wr[f(\a+(\b+\delta_i)t),g(\a+(\b+\delta_i)t)]\leq 0
$$
for all $t\in\RR$.  Taking the limit as $\b\goesto\zero$ and evaluating
at $t=0$ yields
$$
\Wr_i[f,g](\a) = \Wr[f(\a+\delta_i t),g(\a+\delta_i t)]|_{t=0}\leq 0,
$$
by continuity.  Thus (a) implies (c).

To prove that (c) implies (b), let $\a,\b\in\RR^m$ with $\b=(b_1,...,b_m)>\zero$,
and let $a,b\in\RR$ with $b>0$.  By Lemma 2.3, to show that  $g+ y f
\in\StabR[\x,y]$ it suffices to show that
$g(\a+\b t) + (a+\i b) f(\a+\b t) \in\Stab[t]$.  From (c) it follows that
$p=g + a f $ and $q=bf$ are such that $\alpha p +\beta q \in \StabR[\x]$
for all $\alpha,\beta\in\RR$.  By Theorem 2.9, either $p\ll q $ or $q\ll p$.
Now
\showon
\Wr[q(\a+\b t),p(\a+\b t)]
&=&
b \Wr[f(\a+\b t),g(\a+\b t)] \\
&=&
b\sum_{i=1}^m b_i \Wr_i[f,g](\a+\b t) \leq 0,
\showoff
by the Wronskian condition in part (c).  Thus $q(\a+\b t)\ll p(\a +\b t)$,
so that $p(\a+\b t)+\i q(\a +\b t)\in\Stab[t]$.  Since $p+\i q = g + (a+\i b)
f$, this shows that (c) implies (b).
\end{proof}

\begin{exercise}[Corollary 1.10 of \cite{BB5}]
Let $f,g\in\StabR[\x]$ be real stable.  Then $f \ll g$ and $g\ll f$
if and only if $cf=dg$ for some $c,d\in\RR$ not both zero.
\end{exercise}

\begin{prop}[Lemma 3.2 of \cite{BB5}]
Let $V$ be a $\KK$-vector subspace of $\KK[\x]$, with either $\KK=\RR$
or $\KK=\CC$.\\
\textup{(a)}\ If $\KK=\RR$ and $V\subseteq\StabR[\x]$ then $\dim_\RR V \leq
2$.\\
\textup{(b)}\ If $\KK=\CC$ and $V\subseteq\Stab[\x]$ then $\dim_\CC V \leq 1$.
\end{prop}
\begin{proof}
For part (a), suppose to the contrary that $f,g,h\in V$ are linearly
independent over $\RR$ (and hence not identically zero).  By Theorem 2.9,
either $f\ll g$ or $g\ll f$, and similarly for the other pairs $\{f,h\}$ and
$\{g,h\}$.  Renaming these polynomials as necessary, we may assume that
$f\ll h$ and $h\ll g$.  Now, for all $\lambda\in[0,1]$ let $p_\lambda =
(1-\lambda)f+\lambda g$, and note that each $p_{\lambda}\not\equiv 0$.
By Theorem 2.9, for each $\lambda\in[0,1]$ either
$h\ll p_\lambda$ or $p_\lambda \ll h$.  Since $p_0=f\ll h$
and $h\ll g=p_1$, by continuity of the roots of $\{p_\lambda:\
\lambda\in[0,1]\}$ there is a $\lambda\in[0,1]$ such that $h\ll p_\lambda$
and $p_\lambda\ll h$.  But then, by Exercise 2.11, either $\{f,g\}$ is
linearly dependent or $h$ is in the span of $\{f,g\}$, contradicting the
supposition.

For part (b), let $\Re(V)=\{\Re(h):\ h\in V\}$.  Then $\Re(V)$ is a
real subspace of $\StabR[\x]$, so that $\dim_\RR \Re(V)\leq 2$ by part (a).
If $\dim_\RR \Re(V) \leq 1$ then $\dim_\CC V \leq 1$.  In the remaining
case let $\{p,q\}$ be a basis of $\Re(V)$ with $f=p+\i q\in V$.
By Corollary 2.10,  $\Wr_i[q,p](\a)\leq 0$ for all $i\in[m]$ and $\a\in\RR^m$.
Since $p$ and $q$ are not linearly dependent, there is an index $k\in[m]$
such that $\Wr_k[q,p]\not\equiv 0$.

Consider any $g\in V$.  There are reals $a,b,c,d\in\RR$ such that
$$
g = (ap+bq) + \i (cp + dq).
$$
Since $g$ is stable,
$\Wr_k[cp+dq,ap+bq](\a) = (ad-bc)\Wr_k[q,p](\a)\leq 0$
for all $\a\in\RR^m$.  Since $\Wr_k[q,p]\not\equiv 0$, it follows that
$ad-bc\geq 0$.  Now, for any $v,w\in\RR$, $g+(v+\i w)f$ is in $V$.
Since 
$$
g+(v+\i w)f = (a+v)p+(b-w)q +\i((c+w)p+(d+v)q),
$$
this argument shows that $H = (a+v)(d+v) - (b-w)(c+w)\geq 0$
for all $v,w\in\RR$.  But
$$
4H = (2v+a+d)^2 + (2w+c-b)^2 - (a-d)^2 - (b+c)^2,
$$
so that $H\geq 0$ for all $v,w\in\RR$ if and only if $a=d$ and $b=-c$.
This implies that $g=(a+\i c)f$, so that $\dim_\CC V = 1$.
\end{proof}

\section{Multiaffine stable polynomials.}

A polynomial $f$ is \emph{multiaffine} if each indeterminate occurs
at most to the first power in $f$.  For a set $\EuScript{S}$ of polynomials,
let $\EuScript{S}^\ma$ denote the set of multiaffine polynomials in
$\EuScript{S}$.  For multiaffine $f\in\CC[\x]^{\ma}$ and $i\in[m]$
we use the ``ultra-shorthand'' notation $f = f^{i} + x_{i} f_{i}$
in which $f^{i} = f|_{x_{i}=0}$ and $f_{i} = \partial_{i} f$.
This notation is extended to multiple distinct indices in the
obvious way -- in particular,
$$
f = f^{ij} + x_{i} f_{i}^{j} + x_{j} f_{j}^{i} + x_{i} x_{j} f_{ij}.
$$

\subsection{A criterion for real stability.}

For $f\in\CC[\x]$ and $\{i,j\}\subseteq[m]$, let
$$
\Delta_{ij}f
= \partial_i f \cdot \partial_j f - f \cdot \partial_i\partial_j f.
$$
Notice that for $f\in\CC[\x]^\ma$,
$$
\Delta_{ij}f
=
f_{i}^{j} f_{j} - f^{j} f_{ij}
=
\Wr_{i}[f^{j},f_{j}] = -\Wr_{i}[f_{j},f^{j}],
$$
and
$$
\Delta_{ij}f = f_{i}^{j} f_{j}^{i} - f^{ij} f_{ij}.
$$

\begin{theorem}[Theorem 5.6 of \cite{B} and Theorem 3 of \cite{WW}]
Let $f\in\RR[\x]^\ma$ be multiaffine.  The following are equivalent:\\
\textup{(a)}\ $f$ is real stable.\\
\textup{(b)}\ For all $\{i,j\}\subseteq[m]$ and all $\a\in\RR^m$,
$\Delta_{ij}f(\a)\geq 0$.\\
\textup{(c)}\ Either $m=1$, or there exists $\{i,j\}\subseteq[m]$ such that
$f_{i}$, $f^{i}$, $f_{j}$ and $f^{j}$ are real stable, and
$\Delta_{ij}f(\a)\geq 0$ for all $\a\in\RR^m$.
\end{theorem}
\begin{proof}
To see that (a) implies (b), fix $\{i,j\}\subseteq[m]$.
Proposition 2.7(a) shows that $f_{j}\ll f^{j}$, and from the
calculation above and Corollary 2.10, it follows that
$\Delta_{ij}f(\a)  = -\Wr_{i}[f_{j},f^{j}](\a)\geq 0$
for all $\a\in\RR^m$.

We show that (b) implies (a) by induction on $m$, the base
case $m=1$ being trivial.  For the induction step let $f$
be as in part (b), let $a\in\RR$, and let $g=f|_{x_{m}=a}$.
For all $\{i,j\}\subseteq[m-1]$ and $\a\in\RR^{m-1}$,
$\Delta_{ij}g(\a) = \Delta_{ij}f(\a,a)\geq 0$.  By induction,
$g = f^{m} + a f_{m}$ is real stable for all $a\in\RR$; it
follows that $a f_{m} + b f^{m}\in\StabR[x_{1},\ldots,x_{m-1}]$
for all $a,b\in\RR$.  Furthermore, for all $j\in[m-1]$ and
$\a\in\RR^{m-1}$,
$\Wr_{j}[f_{m},f^{m}](\a) = -\Delta_{jm} f(\a,1)\leq 0$.
This verifies condition (c) of Corollary 2.10 for the pair
$(f_{m},f^{m})$, and it follows that $f= f^{m}+ x_{m} f_{m}
\in\StabR[\x]$, completing the induction.

It is clear that (a) and (b) imply (c) -- we show that (c) implies (b) below.
This is clear if $m\leq 2$,  so assume that $m\geq 3$.  To begin with, let
$\{h,i,j\}\subseteq[m]$ be three distinct indices, and consider $\Delta_{ij}f$
as a polynomial in $x_{h}$.  That is, $\Delta_{ij}f =
A_{hij} x_{h}^{2} + B_{hij} x_{h} + C_{hij}$ in which
\showon
A_{hij} &=& f_{hi}^{j} f_{hj}^{i} - f_{h}^{ij} f_{hij} = \Delta_{ij}f_{h},\\
B_{hij} &=& f_{hi}^{j} f_{j}^{hi} - f_{h}^{ij} f_{ij}^{h}
+ f_{i}^{hj} f_{hj}^{i} - f^{hij} f_{hij},\ \mathrm{and}\\
C_{hij} &=& f_{i}^{hj} f_{j}^{hi} - f^{hij} f_{ij}^{h} = \Delta_{ij} f^{h}.
\showoff
If $\Delta_{ij}f(\a)\geq 0$ for all $\a\in\RR^{m}$ then this
quadratic polynomial in $x_{h}$:
$$
\Delta_{ij}f(a_{1},\ldots,a_{h-1},x_{h},a_{h+1},\ldots,a_{m})
$$
has a nonpositive discriminant for all $\a\in\RR^m$.  That is,
$D_{hij} = B_{hij}^{2}-4A_{hij}C_{hij}$
is such that $D_{hij}(\a)\leq 0$ for all $\a\in\RR^{m}$.

It is a surprising fact that as a polynomial in
$\{x_{k}:\ k\in[m]\drop\{h,i,j\}\}$, $D_{hij}$ is invariant under all
six permutations of its indices, as is seen by direct calculation:
\showon
D_{hij}
&=&
(f_{h}^{ij}f_{ij}^{h})^{2} + 
(f_{i}^{hj}f_{hj}^{i})^{2} +
(f_{j}^{hi}f_{hi}^{j})^{2} +
(f_{hij}f^{hij})^{2}\\
& &
-2( f_{h}^{ij}f_{ij}^{h}f_{i}^{hj}f_{hj}^{i}
+ f_{i}^{hj}f_{hj}^{i}f_{j}^{hi}h_{hi}^{j}
+ f_{j}^{hi}f_{hi}^{j}f_{h}^{ij}f_{ij}^{h} )\\
& &
-2 (f_{h}^{ij}f_{ij}^{h} + f_{i}^{hj}f_{hj}^{i} +
f_{j}^{hi} f_{hi}^{j}) f^{hij}f_{hij}\\
& & 
+ 4 f_{ij}^{h} f_{hj}^{i} f_{hi}^{j} f^{hij} 
+ 4 f_{h}^{ij} f_{i}^{hj} f_{j}^{hi} f_{hij}.
\showoff

Now for the proof that (c) implies (b) when $m\geq 3$.
Consider any $h\in[m]\drop\{i,j\}$.  Then
$$
\Delta_{hi}f = A_{jhi} x_{j}^{2} + B_{jhi} x_{j} + C_{jhi}
$$
has discriminant $D_{jhi}=D_{hij}$.  Since $f_{j}$ and
$f^{j}$ are real stable, we have
$A_{jhi}(\a)=\Delta_{hi}f_{j}(\a)\geq 0$ and
$C_{jhi}(\a)=\Delta_{hi}f^{j}(\a)\geq 0$ for all $\a\in\RR^{m}$.
Since $\Delta_{ij} f(\a)\geq 0$ for all $\a\in\RR^{m}$
it follows that $D_{jhi}(\a)=D_{hij}(\a)\leq 0$ for all
$\a\in\RR^{m}$.  It follows that $\Delta_{hi}f(\a)\geq 0$
for all $\a\in\RR^{m}$.  (Note that if $B^{2}-4AC\leq 0$
and either $A=0$ or $C=0$, then $B=0$.)  A similar
argument using the fact that $f_{i}$ and $f^{i}$ are real
stable shows that $\Delta_{hj}f(\a)\geq 0$ for all $\a\in\RR^{m}$.

It remains to show that $\Delta_{hk}f(\a)\geq 0$ for all
$\a\in\RR^{m}$ when $\{h,k\}$ is disjoint from $\{i,j\}$.
We have seen that $\Delta_{hi}f(\a)\geq 0$ for all $\a\in\RR^{m}$,
and we know that both $f_{i}$ and $f^{i}$ are real stable.
The argument above applies once more:\ $\Delta_{hi}f(\a)\geq 0$
for all $\a\in\RR^{m}$, so that $D_{ihk}(\a)=D_{khi}(\a)\leq 0$
for all $\a\in\RR^{m}$, and then since $A_{ihk}(\a)\geq 0$
and $C_{ihk}(\a)\geq 0$ for all $\a\in\RR^{m}$ it follows that
$\Delta_{hk} f(\a)\geq 0$ for all $\a\in\RR^{m}$.  Thus
(c) implies (b).
\end{proof}
\noindent

\subsection{Linear transformations preserving stability -- 
multiaffine case.}

\begin{lemma}[Lieb-Sokal Lemma, Lemma 2.1 of \cite{BB5}]
Let $g(\x)+yf(\x)\in\Stab[\x,y]$ be stable and such that
$\deg_i(f)\leq 1$.  Then $g-\partial_{i}f \in\Stab[\x]$ is stable.
\end{lemma}
\begin{proof}
Since $g$ is stable (by specialization to $y=0$), there is nothing to
prove if $\partial_i f\equiv 0$ identically, so assume otherwise
(and hence that $f\not\equiv 0$).
By permutation we can assume that $i=1$.  Since $f$ is stable and
$z_1,z\in\HH$ imply that $z_1-z^{-1}\in\HH$, it follows that
$$
yf(x_1-y^{-1},x_2,...,x_m)= - \partial_1 f(\x) + yf(\x)
$$
is stable.  Proposition 2.7(b) implies that for all $\z\in\HH^m$,
$$
\Im\left(\frac{g(\z)-\partial_1 f(\z)}{f(\z)}\right) =
\Im\left(\frac{g(\z)}{f(\z)}\right) +
\Im\left(\frac{-\partial_1 f(\z)}{f(\z)}\right) \geq 0.
$$
Thus, by Proposition 2.7(b) again, $g-\partial_1 f + yf$ is stable.
Specializing to $y=0$ shows that $g-\partial_1 f$ is stable.
\end{proof}

\begin{exercise}[Lemma 3.1 of \cite{BB5}]
Let $f\in\CC[\x]^\ma$ and $\w\in\HH^m$.  Then for all $\epsilon>0$
sufficiently small, $(\x+\w)^{[m]}+\epsilon f(\x)$ is stable.
(Here $(\x+\w)^{[m]}=\prod_{i=1}^m(x_i+w_i)$.)
\end{exercise}

For a linear transformation $T:\CC[\x]^{\ma}\goesto\CC[\x]$ of
multiaffine polynomials, define the \emph{algebraic symbol} of $T$ to
be the polynomial
$$
T((\x+\y)^{[m]})
= T\left(\prod_{i=1}^m (x_i+y_i)\right)
= \sum_{S\subseteq[m]} T(\x^{S})\y^{[m]\drop S}
$$
in $\CC[x_{1},\ldots,x_{m},y_{1},\ldots,y_{m}]=\CC[\x,\y]$.

\begin{theorem}[Theorem 1.1 of \cite{BB5}]
Let $T:\CC[\x]^{\ma}\goesto\CC[\x]$ be a linear transformation.
Then $T$ maps $\Stab[\x]^{\ma}$ into $\Stab[\x]$ if and only if 
either\\
\textup{(a)}\ $T(f) = \eta(f)\cdot p$ for some linear functional
$\eta:\CC[\x]^{\ma}\goesto\CC$ and $p\in\Stab[\x]$, or\\
\textup{(b)}\ the polynomial $T((\x+\y)^{[m]})$ is stable in $\Stab[\x,\y]$.
\end{theorem}
\begin{proof}
First, assume (b) that $T((\x+\y)^{[m]})\in\Stab[\x,\y]$ is stable.
By inversion, it follows that $\y^{[m]}T((\x-\y^{-\one})^{[m]})$ is also
stable. Thus, if $f\in\Stab[w_1,...,w_m]$ is stable then
$$
\y^{[m]}T((\x-\y^{-\one})^{[m]})f(\w)
=
\sum_{S\subseteq[m]} T(\x^{S})(-\y)^S f(\w)
$$
is stable.  If $f$ is also multiaffine then repeated application of
the Lieb-Sokal Lemma 3.2 (replacing $y_i$ by $-\partial/\partial w_i$ for
$i\in[m]$) shows that
$$
\sum_{S\subseteq[m]} T(\x^{S})\frac{\del^S}{\del \w^S} f(\w)
$$
is stable.  Finally, specializing to $\w=\zero$ shows that $T(f(\x))$
is stable.  Thus, the linear transformation $T$ maps $\Stab[\x]^\ma$
into $\Stab[\x]$.  This is clearly also the case if (a) holds.

Conversely, assume that $T$ maps $\Stab[\x]^\ma$ into $\Stab[\x]$.
Then for any $\w\in\HH^m$, $(\x+\w)^{[m]}\in\Stab[\x]^\ma$, so that
$T((\x+\w)^{[m]})\in\Stab[\x]$.

First, assume that there is a $\w\in\HH^m$ for which
$T((\x+\w)^{[m]})\equiv 0$ identically.  For any $f\in\CC[\x]^\ma$ let $\epsilon>0$
be as in Exercise 3.3.  Then $\epsilon T(f) = T((\x+\w)^{[m]}+\epsilon f)$
is stable, so that $T(f)$ is stable.  Thus, the image of $\CC[\x]^\ma$ under
$T$ is a $\CC$-subspace of $\Stab[\x]$.  By Proposition 2.12(b), $T$ has the
form of case (a).

Secondly, if $T((\x+\w)^{[m]})\not\equiv 0$ for all $\w\in\HH^m$ then,
since each of these polynomials is in $\Stab[\x]$, we have
$T((\x+\w)^{[m]})|_{\x=\z}\neq 0$ for all $\z\in\HH^m$ and $\w\in\HH^m$.
This shows that $T((\x+\y)^{[m]}))$ is stable in $\Stab[\x,\y]$, which is
the form of case (b).
\end{proof}

Theorem 3.4 has a corresponding real form -- the proof is completely 
analogous.
\begin{theorem}[Theorem 1.2 of \cite{BB5}]
Let $T:\RR[\x]^{\ma}\goesto\RR[\x]$ be a linear transformation.
Then $T$ maps $\StabR[\x]^{\ma}$ into $\StabR[\x]$ if and only if 
either\\
\textup{(a)}\ $T(f) = \eta(f)\cdot p + \xi(f)\cdot q$ for some linear functionals
$\eta,\xi:\RR[\x]^{\ma}\goesto\RR$ and $p,q\in\StabR[\x]$ such that
$p\ll q$, or\\
\textup{(b)}\ the polynomial $T((\x+\y)^{[m]})$ is real stable
in $\StabR[\x,\y]$, or\\
\textup{(c)}\ the polynomial $T((\x-\y)^{[m]})$ is real stable
in $\StabR[\x,\y]$.
\end{theorem}
\begin{proof}
\textbf{Exercise 3.6.}
\end{proof}

\section{The Grace-Walsh-Szeg\H{o} Coincidence Theorem.}

Let $f\in\CC[x]$ be a univariate polynomial of degree at most $m$,
and let $\x=(x_1,...,x_m)$ as usual.  For $0\leq j\leq m$, the
\emph{$j$-th elementary symmetric function} of $\x$ is
$$
e_j(\x) = \sum_{1\leq i_1<\cdots< i_j\leq m} x_{i_1}\cdots x_{i_j}
= \sum_{S\subseteq[m]:\ |S|=j} \x^{S}.
$$
The \emph{$m$-th polarization} of $f$ is the polynomial obtained
as the image of $f$ under the linear transformation $\Pol_m$
defined by $x^j \mapsto \binom{m}{j}^{-1}e_j(\x)$ for all $0\leq j\leq m$,
and linear extension.  In other words, $\Pol_m f$ is the unique
multiaffine polynomial in $\CC[\x]^\ma$ that is invariant under all
permutations of $[m]$ and such that $\Pol_m f(x,...,x) = f(x)$.
A \emph{circular region} is a nonempty subset $\A$ of $\CC$ that is either
open or closed, and which is bounded by either a circle or a straight line.

\begin{theorem}[Grace-Walsh-Szeg\H{o}, Theorem 3.4.1b of \cite{RS}]
Let $f\in\CC[x]$ have degree at most $m$ and let $\A$ be a circular region.
If either $\deg(f)=m$ or $\A$ is convex, then for every $\z\in\A^m$ there
exists $z\in\A$ such that $\Pol_m f(\z)=f(z)$.
\end{theorem}

Figure 1 illustrates the Grace-Walsh-Szeg\H{o} (GWS) Theorem for the polynomial
$f(x)=x^5+10x^2+1$.  The black dots mark the solutions to $f(x)=0$.
Any permutation of the red (grey) dots is a solution to
$\Pol_5 f(x_1,...,x_5)=0$.  By GWS, any circular region containing all the
red dots must contain at least one of the black dots.  The figure indicates
the boundaries of several circular regions for which this condition is met.
\pix{0.8}{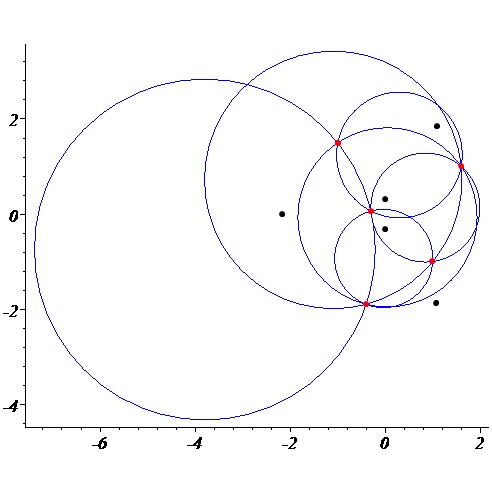}{Illustration of the Grace-Walsh-Szeg\H{o} 
Theorem.}

The proof of GWS in this section is adapted from Borcea and 
Br\"and\'en \cite{BB6}.

\subsection{Reduction to the case of stable polynomials.}

First of all, it suffices to prove GWS for open circular regions, since a closed
circular region is the intersection of all the open circular regions which
contain it.  Second, it suffices to show that for any $g\in\CC[x]$ of degree
at most $m$, if $\deg(g)=m$ or $\A$ is convex, and $\z\in\A^m$ is such that
$\Pol_m g(\z)= 0$, then there exists $z\in\A$ such that $g(z)=0$.
This implies the stated form of GWS by applying this special case to
$g(x)=f(x)-c$, where $c=\Pol_m f(\z)$.  Stated otherwise, it suffices to show
that if $f(z)\neq 0$ for all $z\in\A$ then $\Pol_m f(\z)\neq 0$ for
all $\z\in\A^m$ (provided that either $\deg(f)=m$ or $\A$ is convex).

Let $\MM$ be the set of M\"obius transformations
$z\mapsto\phi(z)=(az+b)/(cz+d)$ with $a,b,c,d\in\CC$ and $ab-cd=\pm 1$.
Then $\MM$ with the operation of functional composition is a group of
conformal transformations of the Riemann sphere $\widehat{\CC}=
\CC\cup\{\infty\}$, and it acts simply transitively on the set of all ordered
triples of distinct points of $\widehat{\CC}$. 
Consequently, for any open circular region $\A$ there is a
$\phi\in\MM$ such that $\phi(\HH)=\{\phi(z):\ z\in\HH\}=\A$.
We henceforth regard circular regions as subsets of $\widehat{\CC}$.
Note that an open circular region $\A$ is convex if and only if
it does not contain $\infty$.  (The point $\infty$ is on the boundary
of any open half-plane.)  In this case,
if $\phi(z)=(az+b)/(cz+d)$ is such that $\phi(\HH)=\A$ then $cz+d\neq 0$
for all $z\in\HH$.

Given $0\not\equiv f\in\CC[x]$ of degree at most $m$, consider the
polynomial $\widetilde{f}(x)=(cx+d)^m f((ax+b)/(cx+d))$.
If either $\deg(f)=m$ or $\A$ is convex, then $f$ is nonvanishing on $\A$
if and only if $\widetilde{f}(z)$ is nonvanishing on $\HH$.
Also,
$$
\Pol_m \widetilde{f}(\x) = 
\Pol_m f (\phi(x_1),...,\phi(x_m)) \cdot \prod_{i=1}^m(cx_i+d).
$$
Thus, to prove GWS it suffices to prove the following lemma.
\begin{lemma}
Let $f\in\CC[x]$ be a univariate polynomial of degree at most
$m$.  Then $\Pol_m f$ is stable if and only if $f$ is stable.
\end{lemma}
\noindent
Clearly, diagonalization implies that if $\Pol_m f$ is stable then $f$
is stable, so only the converse implication needs proof.  This is
accomplished in the following two easy steps.

\subsection{Partial symmetrization.}

The group $\S(m)$ of all permutations $\sigma:[m]\goesto[m]$ acts on $\CC[\x]$
by the rule $\sigma(f)(x_1,...,x_m)=f(x_{\sigma(1)},...,x_{\sigma(m)})$.
Notice that
$$
\sigma(\x^\alpha) = \prod_{i=1}^m x_{\sigma(i)}^{\alpha(i)} =
\prod_{i=1}^m x_i^{\alpha\circ\sigma^{-1}(i)} = \x^{\alpha\circ\sigma^{-1}}.
$$
For $\{i,j\}\subseteq[m]$, let $\tau_{ij}$ be
the transposition that exchanges $i$ and $j$ and fixes all other elements of $[m]$.
\begin{lemma}
Let $0\leq \lambda\leq 1$ and $\{i,j\}\subseteq[m]$, and let
$T_{ij}^{(\lambda)}=(1-\lambda)+\lambda\tau_{ij}$.  If $f\in\Stab[\x]^\ma$ is stable
and multiaffine then $T_{ij}^{(\lambda)} f\in\Stab[\x]^\ma$ is stable and
multiaffine.
\end{lemma}
\begin{proof}
If $f$ is multiaffine then $(1-\lambda)f + \lambda\tau_{ij}(f)$
is also multiaffine.  We apply Theorem 3.4 to show that $T=T_{ij}^{(\lambda)}$
preserves stability of multiaffine polynomials.  By permutation we can assume
that $\{i,j\}=\{1,2\}$. The algebraic symbol of $T$ is
$$
T((\x+\y)^{[m]}) = T((x_1+y_1)(x_2+y_2))\cdot
\prod_{i=3}^m(x_i+y_i).
$$
Clearly, this is stable if and only if the same is true of
$T((x_1+y_1)(x_2+y_2))$.
Exercise 4.4 completes the proof.
\end{proof}
\begin{exercise}
Use the results of Sections 2.4 or 3.1 to show that for
$0\leq \lambda \leq 1$, the polynomial
$$
x_{1}x_{2} + ((1-\lambda)x_{1} + \lambda x_{2})y_{2} +
(\lambda x_{1} + (1-\lambda)x_{2})y_{1} + y_{1}y_{2}
$$
is real stable.
\end{exercise}

\subsection{Convergence to the polarization.}

Let $0\not\equiv f(x)\in\Stab[x]$ be a univariate stable polynomial of degree
at most $m$:\ say $f(x) = c(x-\xi_1)\cdots(x-\xi_n)$ in which $c\neq 0$,
$n\leq m$, and $\xi_i\not\in\HH$ for all $i\in[n]$.  Then the polynomial
$F_0\in\CC[\x]$ defined by
$$
F_0(x_1,...,x_m)=c(x_1-\xi_1)\cdots(x_n-\xi_n)
$$
is multiaffine and stable, and $F_0(x,...,x)=f(x)$.
Let $\Sigma=(\{i_k,j_k\}:\ k\in\NN)$ be a sequence of two-element subsets of
$[m]$, and for each $k\in\NN$ let $T_k = T_{i_k j_k}^{(1/2)}$ and define
$F_{k+1} = T_k(F_k)$.  By induction using Lemma 4.3, each
$F_k\in\Stab[\x]^\ma$ is multiaffine and stable, and
$F_k(x,...,x)=f(x)$ for all $k\in\NN$.  We will construct such
a sequence $\Sigma$ for which $(F_k:\ k\in\NN)$ converges to $\Pol_m f$.

Let $P\in\CC[\x]^\ma$ be multiaffine, say $P(\x)=\sum_{S\subseteq[m]} c(S)\x^S$.
For $\{i,j\}\subseteq [m]$ let
$$
\omega_{ij}(P) = \sum_{S\subseteq[m]} |c(S)-c(\tau_{ij}(S))|
$$
be the \emph{$ij$-th imbalance} of $P$, and let
$||P||=\sum_{\{i,j\}\subseteq[m]} \omega_{ij}(P)$
be the \emph{total imbalance} of $P$.

\begin{exercise}
\textup{(a)}\
Let $(P_k:\ k\in\NN)$ be polynomials in $\CC[\x]^\ma$  for which
there is a $p\in\CC[x]$ such that $P_k(x,...,x)=p(x)$ for all $k\in\NN$.
If $||P_k||\goesto 0$ as $k\goesto 0$, then $(P_k:\ k\in\NN)$ converges to
a limit $P\in\CC[\x]^\ma$, and $||P||=0$.\\
\textup{(b)}\
For $P\in\CC[\x]^\ma$, $||P||=0$ if and only if $P$ is invariant under
all permutations of $[m]$.  Thus, in part (a) the limit is $P=\Pol_m p$.
\end{exercise}

\begin{exercise}
Let $P\in\CC[\x]^\ma$, let $\{i,j\}\subseteq [m]$, and let $Q=T_{ij}^{(1/2)}
P$.\\
\textup{(a)}\ Then $\omega_{ij}(Q)=0$.\\
\textup{(b)}\ If $h\in[m]\drop\{i,j\}$ then
$\omega_{hi}(Q)\leq (\omega_{hi}(P)+\omega_{hj}(P))/2$,
and similarly for $\omega_{hj}(Q)$.\\
\textup{(c)}\ If $\{h,k\}\subseteq[m]\drop\{i,j\}$ then
$\omega_{hk}(Q)=\omega_{hk}(P)$.\\
\textup{(d)}\ Consequently, $||Q|| \leq ||P|| - \omega_{ij}(P).$
\end{exercise}

Now we choose the sequence $\Sigma=(\{i_k,j_k\}:\ k\in\NN)$ as follows:\
for each $k\in\NN$, $\{i_k,j_k\}\subseteq[m]$ is any pair of indices
$\{i,j\}$ for which $\omega_{ij}(F_k)$ attains its maximum value.  Then
$\omega_{i_k j_k}(F_k)\geq\binom{m}{2}^{-1}||F_{k}||$, so that by
Exercise 4.6(d) and induction on $k\in\NN$,
$$
||F_{k+1}||\leq \left(1-\binom{m}{2}^{-1}\right)||F_k||
\leq \left(1-\binom{m}{2}^{-1}\right)^{k+1} ||F_0||.
$$
Thus, by Exercise 4.5, $F_k$ converges to $\Pol_m f$,
the $m$-th polarization of $f$.  Finally, since each $F_k$ is stable 
(and the limit is a polynomial),  Hurwitz's Theorem implies that $\Pol_m f$
is stable.  This completes the proof of Lemma 4.2, and hence of 
Theorem 4.1.

\section{Polarization arguments and stability preservers.}

For $\kappa\in\NN^m$ and a set $\S\subseteq\CC[\x]$ of polynomials,
let $\S^{\leq\kappa}$ be the set of all $f\in\S$ such that $\deg_i(f)\leq
\kappa(i)$ for all $i\in[m]$.  Let
$$
I(\kappa) = \{(i,j):\ i\in[m]\ \mathrm{and}\ j\in[\kappa(i)]\}
$$
and let $\u=\{u_{ij}:\ (i,j)\in I(\kappa)\}$ be indeterminates.
For $f\in\CC[\x]^{\leq\kappa}$, Let $\Pol_{\kappa(i)}^{(i)} f$ denote
the \emph{$\kappa(i)$-th polarization of $x_i$ in $f$}:\ this is the image
of $f$ under the linear transformation $\Pol_{\kappa(i)}^{(i)}$ defined by
$x_i^j \mapsto \binom{\kappa(i)}{j}^{-1} e_j(u_{i1},...,u_{i\kappa(i)})$
for each $0\leq j \leq\kappa(i)$, and linear extension.  Finally, the 
\emph{$\kappa$-th polarization of $f$} is
$$
\Pol_\kappa f = \Pol_{\kappa(m)}^{(m)} \circ\cdots\circ\Pol_{\kappa(1)}^{(1)} f.
$$
This defines a linear transformation $\Pol_{\kappa}:\CC[\x]^{\leq\kappa}
\goesto\CC[\u]^{\ma}$.

\subsection{The real stability criterion revisited.}

\begin{prop}
Let $\kappa\in\NN^{m}$ and $f\in\CC[x]^{\leq\kappa}$.  Then
$\Pol_\kappa f$ is stable if and only if $f$ is stable.
\end{prop}
\begin{proof}
Diagonalization implies that if $\Pol_\kappa f$ is stable then $f$
is stable, so only the converse implication needs proof.  Assume that
$f$ is stable, and let $z_{ij}\in\HH$ for $(i,j)\in I(\kappa)$.
By induction on $m$, repeated application of GWS 
shows that there are $\z=(z_{1},\ldots,z_{m})\in\HH^m$ such that
$$
\Pol_{\kappa} f(z_{ij}:\ (i,j)\in I(\kappa)) = f(\z).
$$
Since $f$ is stable it follows that $\Pol_{\kappa} f$ is stable.
\end{proof}

If $f\in\RR[\x]^{\leq\kappa}$ then Theorem 3.1 applies to $\Pol_\kappa f$.
Thus, Proposition 5.1 bootstraps the real stability criterion from
multiaffine to arbitrary polynomials.  This is a typical application of
the GWS Theorem.

\subsection{Linear transformations preserving stability --
polynomial case.}

\begin{theorem}[Theorem 1.1 of \cite{BB5}]
Let $\kappa\in\NN^m$, and let $T:\CC[\x]^{\leq\kappa}\goesto\CC[\x]$
be a linear transformation.
Then $T$ maps $\Stab[\x]^{\leq\kappa}$ into $\Stab[\x]$ if and only if 
either\\
\textup{(a)}\ $T(f) = \eta(f)\cdot p$ for some linear functional
$\eta:\CC[\x]^{\leq\kappa}\goesto\CC$ and $p\in\Stab[\x]$, or\\
\textup{(b)}\ the polynomial $T((\x+\y)^\kappa)$ is stable in $\Stab[\x,\y]$.
\end{theorem}
\begin{proof}
Let $\u=\{u_{ij}:\ (i,j)\in I(\kappa)\}$, and define a linear 
transformation $\widetilde{T}:\CC[\u]^{\ma}\goesto\CC[\x]$ as follows.
For every $A\subseteq I(\kappa)$, define $\alpha(A):[m]\goesto\NN$ by 
putting $\alpha(A,i) = |\{j\in[\kappa(i)]:\ (i,j)\in A\}|$ for each 
$i\in[m]$.  Then for each $A\subseteq I(\kappa)$ define
$\widetilde{T}(\u^{A}) = T(\x^{\alpha(A)})$, and extend this linearly to
all of $\CC[\u]^{\ma}$.  Let $\Delta:\CC[\u]^{\ma}\goesto\CC[\x]$ be
the diagonalization operator defined by $\Delta(u_{ij})=x_{i}$ for all
$(i,j)\in I(\kappa)$, extended algebraically.

Notice that $T = \widetilde{T}\circ\Pol_{\kappa}$,
and that $\widetilde{T} = T\circ\Delta$.  By Proposition 5.1 (and 
Lemma 2.4), it follows that $T$ preserves stability if and 
only if $\widetilde{T}$ preserves stability.  This is
equivalent to one of two cases in Theorem 3.4.

In case (a), if $\widetilde{T}= p\cdot\widetilde{\eta}$ for some 
$p\in\Stab[\x]$ and linear functional $\widetilde{\eta}:\CC[\y]^{\ma}\goesto\CC$
then $T=p\cdot (\eta\circ \Pol_{\kappa})$ is also in case (a).
Conversely, if $T$ is in case (a) then the same is true of 
$\widetilde{T}$, by construction.

In case (b), let 
$\Pol_{\kappa}^{(\y)}:\CC[\y]^{\leq\kappa}\goesto\CC[\v]^{\ma}$
denote the $\kappa$-th polarization of the $\y$ variables.
The symbols of $T$ and $\widetilde{T}$ are related by
$$
\widetilde{T}((\u+\v)^{I(\kappa)}) =
(T\circ \Delta)((\u+\v)^{I(\kappa)}) =
\Pol_{\kappa}^{(\y)}T((\x+\y)^{\kappa}),
$$
and Proposition 5.1 shows that $T$ is in case (b) if and only if
$\widetilde{T}$ is in case (b).
\end{proof}

\subsection{Linear transformations preserving stability -- 
transcendental case.}

\begin{exercise}
Let $T:\CC[\x]\goesto\CC[\x]$ be a linear transformation.\\
\textup{(a)}\ Then $T:\Stab[\x]\goesto\Stab[\x]$ if and only if
$T:\Stab[\x]^{\leq\kappa}\goesto\Stab[\x]$ for all 
$\kappa\in\NN^{m}$.\\ 
\textup{(b)}\ Define $S:\CC[\x,\y]\goesto\CC[\x,\y]$ by $S(\x^{\alpha}\y^{\beta})=
T(\x^{\alpha})\y^{\beta}$ and linear extension.  If $T((\x+\u)^{\kappa})$
is stable for all $\kappa\in\NN^{m}$ then
$S((\x+\u)^{\kappa}(\y+\v)^{\beta})$ is stable for all $\kappa,\beta\in\NN^{m}$.
\end{exercise}

Let $\overline{\Stab[\x]}$ denote the set of all power series in
$\CC[[\x]]$ that are obtained as the limit of a sequence of stable polynomials
in $\Stab[\x]$ which converges uniformly on compact sets.
Theorem 5.4 is an astounding generalization of the P\'olya-Schur Theorem.
For $\alpha\in\NN^{m}$, let $\alpha!=\prod_{i=1}^{m}\alpha(i)!$.

\begin{theorem}[Theorem 1.3 of \cite{BB5}]
Let $T:\CC[\x]\goesto\CC[\x]$ be a linear transformation.
Then $T$ maps $\Stab[\x]$ into $\Stab[\x]$ if and only if 
either\\
\textup{(a)}\ $T(f) = \eta(f)\cdot p$ for some linear functional
$\eta:\CC[\x]\goesto\CC$ and $p\in\Stab[\x]$, or\\
\textup{(b)}\ the power series
$$T(\e^{-\x\y}) = \sum_{\alpha:[m]\goesto\NN}(-1)^{\alpha}
T(\x^\alpha)\frac{\y^\alpha}{\alpha!}
$$
is in $\overline{\Stab[\x,\y]}$ 
\end{theorem}
(Theorem 3.5 has a similar extension -- see Theorems 1.2 and 1.4 of
\cite{BB5}.)

For $\alpha\leq \beta$ in $\NN^{m}$, let $(\beta)_{\alpha} = 
\beta!/(\beta-\alpha)!$, and for $\alpha\not\leq\beta$ let
$(\beta)_{\alpha} = 0$.

\begin{theorem}[Theorem 5.1 of \cite{BB5}]
Let $F(\x,\y)= \sum_{\alpha\in\NN^{m}} P_{\alpha}(\x)\y^{\alpha}$ be a
power series in $\CC[\x][[\y]]$ (so that each $P_{\alpha}\in\CC[\x]$).
Then $F(\x,\y)$ is in $\overline{\Stab[\x,\y]}$ if and only if
for all $\beta\in\NN^{m},$
$$
\sum_{\alpha\leq\beta}(\beta)_{\alpha}P_{\alpha}(\x)\y^{\alpha}
$$
is stable in $\Stab[\x,\y]$.
\end{theorem}
(This implies the analogous result for real stability, since
$\StabR[\x]=\Stab[\x]\cap\RR[\x]$.)

\begin{exercise}
Derive Theorem 5.4 from Theorems 5.2 and 5.5.  (Hint:\ 
$T((\x+\y)^\kappa)$ is stable if and only if $T((\one-\x\y)^\kappa)$
is stable.)
\end{exercise}

One direction of Theorem 5.5 is relatively straightforward.

\begin{lemma}[Lemma 5.2 of \cite{BB5}]
Fix $\beta\in\NN^{m}$.  The linear transformation $T:\y^{\alpha}\mapsto
(\beta)_{\alpha}\y^{\alpha}$ on $\CC[\y]$ preserves stability.
\end{lemma}
\begin{proof}
By Theorem 5.2 and Exercise 5.3(a), it suffices to show that for
all $\kappa\in\NN^{m}$, the polynomial $T((\y+\u)^{\kappa})$
is stable.  Now
$$
T((\y+\u)^{\kappa}) = \prod_{i=1}^{m}
\left[
\sum_{j=0}^{\kappa(i)} j! \binom{\kappa(i)}{j}\binom{\beta(i)}{j} 
y_{i}^{j}u_{i}^{\kappa(i)-j}
\right],
$$
so it suffices to show that for all $k,b\in\NN$, the polynomial
$f(t) = \sum_{j=0}^{k} j! \binom{k}{j}\binom{b}{j} t^{j}$
is real stable.  Let $g(t) = (1+d/dt)^{k} t^{b}$. 
One can check that $f(t) = t^{b} g(1/t)$.  It thus suffices to show
that $1+d/dt$ preserves stability.  For any $a\in\NN$,
$(1+d/dt)(t+u)^{a} = (t+u+a)(t+u)^{a-1}$ is stable, and so Theorem 5.2
implies the result.
\end{proof}

Now, let $F=F(\x,\y)$ be as in the statement of Theorem 5.5, and let
$(F_{n}:\ n\in\NN)$ be a sequence of stable polynomials
$F_{n}(\x,\y)= \sum_{\alpha\in\NN^{m}} P_{n,\alpha}(\x)\y^{\alpha}$ 
in $\Stab[\x,\y]$ converging to $F$ uniformly on compact sets.
Fix $\beta\in\NN$ and define a linear transformation
$T:\CC[\x,\y] \goesto\CC[\x,\y]$ by $T(\x^{\gamma}\y^{\alpha})
= (\beta)_{\alpha}  \x^{\gamma}\y^{\alpha}$ and linear extension.
By Lemma 5.7 and Exercise 5.3, $T$ preserves stability in 
$\Stab[\x,\y]$.  Thus, $(T(F_{n}):\ n\in\NN)$ is a sequence of stable
polynomials converging to $T(F)$.  Since $T(F)$ is a 
polynomial the convergence is uniform on compact sets, and so
Hurwitz's Theorem implies that $T(F)$ is stable.

The converse direction of Theorem 5.5 is considerably more technical,
although the idea is simple.  With $F$ as in the theorem, for each 
$n\geq 1$ let
$$
F_{n}(\x,\y) = \sum_{\alpha \leq n\one}
(n\one)_{\alpha}P_{\alpha}(\x)\frac{\y^{\alpha}}{n^{\alpha}}.
$$
The sequence $(F_{n}:\ n\geq 1)$ converges to $F$,
since for each $\alpha\in\NN^m$,
$n^{-\alpha}(n\one)_\alpha \goesto 1$ as $n\goesto\infty$.
Each $F_{n}$ is stable, by hypothesis (and scaling).  The hard
work is involved with showing that the convergence is uniform on
compact sets.  To do this, Borcea and Br\"and\'en develop a very
flexible multivariate generalization of the Szasz Principle
\cite[Theorem 5.6]{BB5} -- in itself an impressive accomplishment.
Unfortunately, we have no space here to develop this result -- see
Section 5.2 of \cite{BB5}.

\section{Johnson's Conjectures.}

Let $\A = (A_{1},\ldots,A_{k})$ be a $k$-tuple of $n$-by-$n$ matrices.
Define the \emph{mixed determinant} of $\A$ to be
$$
\Det(\A) = \Det(A_{1},\ldots,A_{k}) = 
\sum_{(S_{1},\ldots,S_{k})} \prod_{i=1}^{k} \det A_{i}[S_{i}],
$$
in which the sum is over all ordered sequences of $k$ pairwise disjoint
subsets of $[n]$ such that $[n]=S_{1}\cup\cdots\cup S_{k}$,  and
$A_{i}[S_{i}]$ is the principal submatrix of $A_{i}$ supported on rows 
and columns in $S_{i}$.  Let $A_i(S_i)$ be the complementary
principal submatrix supported on rows and columns not in $S_i$,
and for $j\in[n]$ let $A_i(j)=A_i(\{j\})$.

For example, when $k=2$ and $A_{1}=xI$ and $A_{2}=-B$, this specializes 
to $\Det(xI,-B)=\det(xI-B)$, the characteristic polynomial of $B$.
In the late 1980s, Johnson made three conjectures about
the $k=2$ case more generally.\\

\noindent
\textbf{Johnson's Conjectures.}\  \emph{
Let $A$ and $B$ be $n$-by-$n$ matrices, with $A$ positive definite
and $B$ Hermitian.\\
\textup{(a)}\  Then $\Det(xA,-B)$ has only real roots.\\
\textup{(b)}\  For $j\in[n]$, the roots of 
$\Det(xA(j),-B(j))$ interlace those of $\Det(xA,-B)$.\\
\textup{(c)}\  The inertia of $\Det(xA,-B)$ is the same as that of
$\det(xI-B)$.}\\

In part (c), the \emph{inertia} of a univariate real stable polynomial $p$
is the triple $\iota(p)=(\iota_{-}(p),\iota_{0}(p),\iota_{+}(p))$
with entries the number of negative, zero, or positive roots of $p$, respectively.

In 2008, Borcea and Br\"and\'en \cite{BB1} proved all three of these
statements in much greater generality.

\begin{theorem}[Theorem 2.6 of \cite{BB1}]
Fix integers $\ell,m,n\geq 1$.  For $h\in[\ell]$ and $i\in[m]$
let $B_{h}$ and $A_{hi}$ be $n$-by-$n$ matrices, and let
$$
L_{h} = \sum_{i=1}^{m} x_{i} A_{hi} + B_{h}.
$$
\textup{(a)}\
If all the $A_{hi}$ are positive semidefinite and all the $B_{h}$
are Hermitian, then $\Det(\L)=\Det(L_{1},\ldots,L_\ell)
\in\StabR[\x]$ is real stable.\\
\textup{(b)}\
For each $j\in[n]$, let $\L(j)=(L_{1}(j),\ldots,L_\ell(j))$.
With the hypotheses of part (a),
the polynomial $\Det(\L)+y\Det(\L(j))\in\StabR[\x,y]$ is real stable.
\end{theorem}
\begin{proof}
Let $Y=\diag(y_1,...,y_n)$ be a diagonal matrix of indeterminates.
By Proposition 2.1, for each $h\in[\ell]$ the polynomial
$$
\det(Y+L_h) = \sum_{S\subseteq[n]}\y^S \det L_h(S)
$$
is real stable in $\StabR[\x,\y]$.
By inversion of all the $\y$ indeterminates, each
$$
\det(I-YL_h)=  \sum_{S\subseteq[n]}(-1)^{|S|} \y^S \det L_h[S]
$$
is real stable.  Since $\prod_{h=1}^\ell \det(I-YL_h)$ is real 
stable, contraction and specialization imply that
$$
\Det(\L) =
(-1)^n 
\left.\frac{\partial^n}{\partial y_1\cdots \partial y_n}
\prod_{h=1}^\ell \det(I-YL_h)\right|_{\y=\zero}
$$
is real stable, proving part (a).

For part (b), let $V$ be the $n$-by-$n$ matrix with all entries
zero except for $V_{jj}=y$.  By part (a),
$$
\Det(V,L_1,...,L_h)=\Det(\L) + y \Det(\L(j))
$$
is real stable.
\end{proof}

Theorem 6.1 (with Corollary 2.10) clearly settles Conjectures (a) and (b).

\begin{proof}[Proof of Conjecture (c).]
Let $A$ and $B$ be $n$-by-$n$ matrices with 
$A$ positive definite and $B$ Hermitian.  Let $(\iota_-,\iota_0,\iota_+)$
be the inertia of $\det(xI-B)$.  Let $f(x) = \Det(xA,-B)$, and let
$(\nu_-,\nu_0,\nu_+)$ be the inertia of $f$.

We begin by showing that $\nu_0=\iota_0$.
Since $\iota_0=\min\{|S|:\ S\subseteq[n]\ \mathrm{and}\ \det(B(S))\neq 0\}$,
it follows that $\nu_0\geq \iota_0$.  The constant term of $f(x)$ is
$(-1)^n\det(B)$, so that if $\iota_0=0$ then $\nu_0=0$.  If $\iota_0=k>0$
then let $S=\{s_1,...,s_k\}\subseteq[n]$ be such that $\det(B(S))\neq 0$.
For $0\leq i\leq k$ let 
$f_i(x)=\Det(A(\{s_1,..,s_i\}),-B(\{s_1,..,s_i\}))$,
so that $f_0(x)=f(x)$.  By Theorem 6.1, the roots of $f_{i-1}$ and of
$f_i$ are interlaced, for each $i\in[k]$.  Thus,
$$
\nu_0 = \iota_0(f_0) \leq \iota_0(f_1)+1 \leq \iota_0(f_2)+2 \leq
\cdots \leq \iota_0(f_k)+k = k = \iota_0,
$$
since $\iota_0(f_k)=0$ because $\det(B(S))\neq 0$.  Therefore $\nu_0=\iota_0$.

For any positive definite matrix $A$, $\Det(xA,-B)$ is a polynomial of degree
$n$.  Suppose that $A$ is such a matrix for which $\nu_+\neq\iota_+$.
Consider the matrices $A_\lambda = (1-\lambda)I+\lambda A$ for
$\lambda\in[0,1]$.  Each of these matrices is positive definite.
From the paragraph above, each of the polynomials
$g_\lambda(x)=\Det(xA_\lambda,-B)$ has $\iota_0(g_\lambda)=\iota_0$.
Since $\iota_+(g_0)=\iota_+\neq\nu_+=\iota_+(g_1)$ and the roots of $g_\lambda$
vary continuously with $\lambda$, there is some value $\mu\in(0,1)$ for which
$\iota_0(g_\mu)>\iota_0$.  This contradiction shows that $\nu_+=\iota_+$,
and hence $\nu_-=\iota_-$ as well.
\end{proof}

Borcea and Br\"and\'en \cite{BB1} proceed to derive many inequalities for
the principal minors of positive semidefinite matrices, and some for merely
Hermitian matrices.  These are applications of inequalities valid more
generally for real stable polynomials.  The simplest of these inequalities are
as follows.

For an $n$-by-$n$ matrix $A$, the \emph{$j$-th symmetrized Fisher product}
is
$$
\sigma_j(A) = \sum_{S\subseteq[n]:\ |S|=j} \det(A[S])\det(A(S)).
$$
and the \emph{$j$-th averaged Fisher product} is
$\widehat{\sigma}_j(A) = \binom{n}{j}^{-1} \sigma_j(A)$.
Notice that $\sigma_j(A)=\sigma_{n-j}(A)$ for all $0\leq j\leq n$.

\begin{coro}
Let $A$ be an $n$-by-$n$ positive semidefinite matrix.\\
\textup{(a)}\ Then $\widehat{\sigma}_j(A)^2 \geq 
\widehat{\sigma}_{j-1}(A) \widehat{\sigma}_{j+1}(A)$
for all $1\leq j\leq n-1$.\\
\textup{(b)}\ Also, $\widehat{\sigma}_0(A)\leq \widehat{\sigma}_1(A)
\leq \cdots \leq \widehat{\sigma}_{\lfloor n/2\rfloor}$.\\
\textup{(c)}\  If $A$ is positive definite and $\det(A)=d$ then
$$
\frac{\widehat{\sigma}_1(A)}{d} \geq
\left(\frac{\widehat{\sigma}_2(A)}{d}\right)^{1/2} \geq
\left(\frac{\widehat{\sigma}_3(A)}{d}\right)^{1/3} \geq \cdots \geq
\left(\frac{\widehat{\sigma}_n(A)}{d}\right)^{1/n} =1.
$$
\end{coro}
\begin{proof}
It suffices to consider positive definite $A$.
By Theorem 6.1, the polynomial $\Det(xA,-A)=\sum_{j=0}^n (-1)^j
\sigma_j(A) x^j$ has only real roots, and these roots are all positive.
Part (a) follows from Newton's Inequalities \cite[Theorem 51]{HLP}.
Part (a) and the symmetry $\sigma_j(A)=\sigma_{n-j}(A)$ for all $0\leq j\leq n$
imply part (b).  Part (c) follows from Maclaurin's Inequalities
\cite[Theorem 52]{HLP}.
\end{proof}

\section{The symmetric exclusion process.}

This section summarizes an application of stable polynomials to probability
and statistical mechanics from a 2009 paper of Borcea, Br\"and\'en and Liggett
\cite{BBL}.

Let $\Lambda$ be a set of \emph{sites}.  A symmetric exclusion process
(SEP) is a type of Markov chain with state space a subset of $\{0,1\}^\Lambda$.
In a state $S:\Lambda\goesto\{0,1\}$, the sites in
$S^{-1}(1)$ are \emph{occupied} and the sites in $S^{-1}(0)$ are \emph{vacant}.
This is meant to model a physical system of particles interacting by means of
hard-core exclusions.  Such models come in many varieties -- to avoid technicalities
we discuss only the case of a finite system $\Lambda$ and continuous time $t$.
(The results of this section extend to countable $\Lambda$ under a reasonable
finiteness condition on the interaction rates.)  Symmetry of the interactions turns 
out to be crucial, but particle number conservation is unimportant.

Let $E$ be a set of two-element subsets of $\Lambda$.  For each $\{i,j\}\in E$,
let $\lambda_{ij}>0$ be a positive real, and let $\tau_{ij}:\Lambda
\goesto\Lambda$ be the permutation that exchanges $i$ and $j$ and fixes all 
other sites.  Our SEP Markov chain $\MM$ proceeds as follows.
Each $\{i,j\}\in E$ has a Poisson process ``clock''
of rate $\lambda_{ij}$, and these are independent of one another.
With probability one, no two clocks ever ring at the same time.  When the clock
of $\{i,j\}$ rings, the current state $S$ is updated to the new state
$S\circ\tau_{ij}$.  In other words, when the $\{i,j\}$ clock rings,
if exactly one of the sites $\{i,j\}$ is occupied then a particle hops
from the occupied to the vacant of these two sites.

Let $\Lambda=[m]$ and $\Omega=\{0,1\}^\Lambda$,
let $\varphi_0$ be an initial probability distribution on $\Omega$,
and let $\varphi_t$ be the distribution of the state of $\MM$, starting
at $\varphi_0$, after evolving for time $t\geq 0$.  We are concerned with
properties of the distribution $\varphi_t$ that hold for all $t\geq 0$.

\subsection{Negative correlation and negative association.}

Consider a probability distribution $\varphi$ on $\Omega$.
An \emph{event} $\E$ is any subset of $\Omega$.  The probability of the event
$\E$ is $\Pr[\E] = \sum_{S\in\E} \varphi(S)$.
An event $\E$ is \emph{increasing} if whenever $S\leq S'$ in $\Omega$
and $S\in\E$, then $S'\in\E$.  For example, if $K$ is any subset of 
$\Lambda$ and $\E_{K}$ is the event that all sites in $K$ are occupied, 
then $\E_{K}$ is an increasing event.  Notice that this event has the 
form $\E_{K}=\E' \times \{0,1\}^{\Lambda\drop K}$ for some event 
$\E'\subseteq\{0,1\}^{K}$.  Two events $\E$ and $\F$ are 
\emph{disjointly supported} when one can partition $\Lambda = A\cup 
B$ with $A\cap B=\none$ and $\E=\E' \times \{0,1\}^{B}$
and $\F=\{0,1\}^{A}\times\F'$ for some events 
$\E'\subseteq\{0,1\}^{A}$ and $\F'\subseteq\{0,1\}^{B}$.

A probability distribution on $\Omega$ is \emph{negatively associated}
(NA) when $\Pr[\E\cap\F]\leq \Pr[\E]\cdot\Pr[\F]$ for any two
increasing events that are disjointly supported.  It is 
\emph{negatively correlated} (NC) when $\Pr[\E_{\{i,j\}}]\leq
\Pr[\E_{\{i\}}]\cdot\Pr[\E_{\{j\}}]$  for any two distinct sites 
$\{i,j\}\subseteq \Lambda$.  Clearly NA implies NC.

It is useful to find conditions under which NC implies NA, since NC
is so much easier to check.  The following originates with Feder and 
Mihail, but many others have contributed their insights -- see Section
4.2 of \cite{BBL}.  The \emph{partition function} of any
$\varphi:\Omega\goesto\RR$ is the real multiaffine polynomial
$$
Z(\varphi)=Z(\varphi;\x) = \sum_{S\in\Omega} \varphi(S) \x^{S}
$$
in $\RR[\x]^{\ma}$.
If $\varphi$ is nonzero and nonnegative, then for any $\a\in\RR^{\Lambda}$
with $\a>\zero$, this defines a probability distribution
$\varphi^{\a}:\Omega\goesto[0,1]$  by setting $\varphi^{\a}(S) =
\varphi(S)\a^{S}/Z(\varphi;\a)$ for all $S\in\Omega$.\\

\noindent
\textbf{Feder-Mihail Theorem} (Theorem 4.8 of \cite{BBL})\textbf{.}\
\emph{
Let $\S$ be a class of nonzero nonnegative functions satisfying the following
conditions.\\
\textup{(i)}\  Each $\varphi\in\S$ has domain $\{0,1\}^{\Lambda}$
for some finite set $\Lambda=\Lambda(\varphi)$.\\
\textup{(ii)}\  For each $\varphi\in\S$, $Z(\varphi)$ is a homogeneous 
polynomial.\\
\textup{(iii)}\  For each $\varphi\in\S$ and $i\in\Lambda(\varphi)$,
$Z(\varphi)|_{x_{i}=0}$ and $\partial_{i}Z(\varphi)$ are partition
functions of members of $\S$.\\
\textup{(iv)}\  For each $\varphi\in\S$ and $\a\in\RR^{\Lambda(\varphi)}$
with $\a>\zero$, $\varphi^{\a}$ is NC.\\
Then for every $\varphi\in\S$ and $\a\in\RR^{\Lambda(\varphi)}$
with $\a>\zero$, $\varphi^{\a}$ is NA.}

\subsection{A conjecture of Liggett and Pemantle.}

In the early 2000s, Liggett and Pemantle arrived independently at
the following conjecture, now a theorem.

\begin{theorem}[Theorem 5.2 of \cite{BBL}]
If the initial distribution $\varphi_0$ of a SEP is deterministic (\emph{i.e.}
concentrated on a single state) then $\varphi_t$ is NA for all $t\geq 0$.
\end{theorem}
\begin{proof}
This amounts to finding a class $\Stab$ of probability distributions
such that:\\
(1)\ deterministic distributions are in $\Stab$,\\
(2)\ being in $\Stab$ implies NA, and\\
(3)\ time evolution of the SEP preserves membership in $\Stab$.

Borcea, Br\"and\'en, and Liggett \cite{BBL} identified such a 
class:\ $\varphi$ is in $\Stab$ if and only if the partition function
$Z(\varphi)$ is homogeneous, multiaffine, and real stable.
(Notice that if $\varphi$ is in $\Stab$ then $\varphi^{\a}$ is in 
$\Stab$ for all $\a\in\RR^{\Lambda}$ with $\a>\zero$, by scaling.)
We proceed to check the three claims above.

Claim (1) is trivial, since if $\varphi(S)=1$ then 
$Z(\varphi)=\x^{S}$, which is clearly homogeneous, multiaffine, and real
stable.

To check claim (2) we verify the hypotheses of the Feder-Mihail 
Theorem.  Hypotheses (i) and (ii) hold since $Z(\varphi)$ is 
multiaffine and homogeneous.  By specialization and contraction,
(iii) holds.  To check (iv), let $\a\in\RR^{\Lambda}$ with $\a>\zero$,
let $\{i,j\}\subseteq \Lambda$, and consider the probability 
distribution $\varphi^{\a}$ on $\Omega$.  The occupation 
probability for site $i$ is
$$
\Pr[\E_{\{i\}}] = \sum_{S\in\{0,1\}^\Lambda:\ S(i)=1}
\frac{\varphi(S)\a^{S}}{Z(\varphi;\a)}
=
a_{i} \frac{\partial_{i}Z(\varphi;\a)}{Z(\varphi;\a)},
$$
and similarly for $\Pr[\E_{\{j\}}]$.  Likewise,
$\Pr[\E_{\{i,j\}}] =
a_{i}a_{j}Z(\varphi;\a)^{-1} \cdot \partial_{i}\partial_{j} Z(\varphi;\a)$.
Now
$$
\Pr[\E_{\{i,j\}}]  - \Pr[\E_{\{i\}}]\cdot \Pr[\E_{\{j\}}] 
=
- \frac{a_{i}a_{j}}{Z(\varphi;\a)^{2}} \cdot \Delta_{ij} Z(\varphi;\a) \leq 0,
$$
by Theorem 3.1.  Thus $\varphi^{\a}$ is NC.
By the Feder-Mihail Theorem, every $\varphi$ in $\Stab$ is NA.

To check claim (3) we need some of the theory of continuous time Markov 
chains.  The time evolution of a Markov chain $\MM$ with finite
state space $\Omega$ is governed by a one-parameter semigroup $T(t)$
of transformations of $\RR^\Omega$.  For a function $F\in\RR^\Omega$
and time $t\geq 0$ and state $S\in \Omega$, $(T(t)F)(S)$ is the
expected value of $F$ at time $t$, given that the initial distribution
of $\MM$ is concentrated at $S$ with probability one at time $0$.
In particular, $\varphi_{t} = T(t)\varphi_{0}$ for all $t\geq 0$, and
all initial distributions $\varphi_{0}$.
In the case of the SEP we are considering, the infinitesimal generator
$\L$ of the semigroup $T(t)$ is given by
$$
\L = \sum_{\{i,j\}\in E} \lambda_{ij} \left( \tau_{ij} - 1 \right).
$$
For each $\{i,j\}\in E$, this replaces each $S\in\Omega$
by $S\circ\tau_{ij}$ at the rate $\lambda_{ij}$.

In preparation for Section 7.3, it is useful to regard $\L$ as an element
of the real semigroup algebra $\Alg=\RR[\EE]$ of the semigroup
$\EE$ of all endofunctions $\ff:\Omega\goesto\Omega$ (with the
operation of functional composition).  The left action of $\EE$ on
$\Omega$ is extended to a left action of $\Alg$ on $\CC[\x]$ as usual:\
for $\ff\in\EE$ and $S\in\Omega$, $\ff(\x^S) = \x^{\ff(S)}$,
extended bilinearly to all of $\Alg$ and $\CC[\x]$.  A permutation
$\sigma\in\S(\Lambda)$ is identified with the endofunction $\ff_\sigma:S
\mapsto S\circ\sigma^{-1}$, so this action of $\Alg$ agrees with
the action of $\S(m)$ in Section 4.2.  A left action of $\Alg$ on $\RR^\Omega$
is defined by $Z(\ff(F))=\ff(Z(F))$ for all $\ff\in\EE$ and $F\in\RR^\Omega$,
and linear extension.  More explicitly, for $\ff\in\EE$, $F\in\RR^\Omega$,
and $S\in\Omega$,
$$
(\ff(F))(S) = F(\ff^{-1}(S)) =
\sum \{F(S'):\ S'\in\Omega\ \ \mathrm{and}\ \ \ff(S')=S\}.
$$

Consider an element of $\Alg$ of the form
$\L = \sum_{i=1}^N \lambda_i (\ff_i-1)$ with all $\lambda_i>0$.
Let $\lambda_i\leq L$ for all $i\in[N]$, and let $K=\sum_{i=1}^N \lambda_i$.
The power series
$$
\exp(t\L) = \e^{-Kt}\sum_{n=0}^{\infty}\frac{t^{n}}{n!}
\left[
\sum_{i=1}^N \lambda_i \ff_i
\right]^{n}
=\sum_{\ff\in\EE} P_{\ff}(t)\cdot\ff
$$
in $\Alg[[t]]$ is such that for each $\ff\in\EE$,
$P_{\ff}(t)\in\RR[[t]]$ is dominated coefficientwise by $\exp((LN-K)t)$.  Thus
$\exp(t\L)\in\Alg[[t]]$ converges for all $t\geq 0$.
The semigroup of transformations generated by $\L$ is $\exp(t\L)$.

To check claim (3) we will show that the semigroup $T(t)$ of the SEP
preserves stability for all $t\geq 0$:\ that is,
if $Z(\varphi_0)$ is stable then $Z(\varphi_t)=T(t)Z(\varphi_0)$
is stable for all $t\geq 0$.
This reduces to the case of a single pair $\{i,j\}\in E$, as follows.
If $\MM_{1}$ and $\MM_{2}$ are Markov chains
on the same finite state space, with semigroups $T_{1}(t)$ and $T_{2}(t)$
generated by $\L_{1}$ and $\L_{2}$, then the semigroup generated by
$\L_{1}+\L_{2}$ is
$$
T(t) = \lim_{n\goesto\infty} \left[ T_{1}(t/n) T_{2}(t/n) \right]^{n},
$$
by the Trotter product formula.  By Hurwitz's Theorem, It follows that
if $T_{i}(t)$ preserves stability for all $t\geq 0$ and $i\in\{1,2\}$,
then $T(t)$ preserves stability for all $t\geq 0$.  By repeated application
of this argument, in order to show that the SEP semigroup $T(t)=\exp(t\L)$
preserves stability for all $t\geq 0$ it is enough to show that
for each $\{i,j\}\in E$, $T_{ij}(t) = \exp(t\lambda_{ij}(\tau_{ij}-1))$
preserves stability for all $t\geq 0$.  Now, since $\tau_{ij}^{2}=1$,
$$
T_{ij}(t)
=
\left(\frac{1+\e^{-2\lambda_{ij}t}}{2}\right)+
\left(\frac{1-\e^{-2\lambda_{ij}t}}{2}\right)\cdot\tau_{ij}.
$$
By Lemma 4.3, this preserves stability for all $t\geq 0$.
This proves Theorem 7.1.
\end{proof}

\subsection{Further observations.}

In verifying the hypotheses of the Feder-Mihail Theorem we used the 
fact that if $f\in\StabR[\x]^{\ma}$ is multiaffine and real stable, 
then $\Delta_{ij}f(\a)\geq 0$ for all $\{i,j\}\subseteq E$ and
$\a\in\RR^{m}$, by Theorem 3.1.  In fact, we only needed the
weaker hypothesis that $\Delta_{ij}f(\a)\geq 0$ for all $\{i,j\}\subseteq E$
and $\a\in\RR^{m}$ with $\a>\zero$.  A multiaffine real polynomial 
satisfying this weaker condition is a \emph{Rayleigh} polynomial.
(This terminology is by analogy with the Rayleigh monotonicity 
property of electrical networks -- see Definition 2.5 of \cite{BBL}
and the references cited there.  Multiaffine real stable polynomials 
are also called \emph{strongly Rayleigh}.)  The class of probability
distributions $\varphi$ such that $Z(\varphi)$ is homogeneous, multiaffine,
and Rayleigh meets all the conditions of the Feder-Mihail Theorem.
It follows that all such distributions are NA.

Claim (2) above can be generalized in another way --
the hypothesis of homogeneity can be removed, as follows.
Let $\y=(y_{1},\ldots,y_{m})$ and let $e_{j}(\y)$
be the $j$-th elementary symmetric function of the $\y$.
Given a multiaffine polynomial $f=\sum_{S\subseteq[m]}c(S)\x^{S}$,
the \emph{symmetric homogenization} of $f$ is the polynomial 
$f_{\mathrm{sh}}(\x,\y)\in\CC[\x,\y]^{\ma}$ defined by
$$
f_{\mathrm{sh}}(\x,\y) = \sum_{S\subseteq[m]} c(S)\x^{S}\binom{m}{|S|}^{-1}
e_{m-|S|}(\y).
$$
Note that $f_{\mathrm{sh}}$ is homogeneous of degree $m$, and
$f_{\mathrm{sh}}(\x,\one)=f(\x)$.

\begin{prop}[Theorem 4.2 of \cite{BBL}]
If $f\in\StabR[\x]^{\ma}$ is multiaffine and real stable then
$f_{\mathrm{sh}}\in\StabR[\x,\y]^{\ma}$ is homogeneous, multiaffine and real stable.
\end{prop}
(We omit the proof.)

\begin{coro}[Theorem 4.9 of \cite{BBL}]
Let $\varphi:\Omega\goesto[0,\infty)$ be such that
$Z(\varphi)$ is nonzero, multiaffine, and real stable.  Then for all 
$\a\in\RR^{m}$ with $\a>0$, $\varphi^{\a}$ is NA.
\end{coro}
\begin{proof}
By Proposition 7.2, $Z_{\mathrm{sh}}(\varphi;\x,\y)$ is nonzero, homogeneous,
multiaffine, and real stable.  This is the partition function for
$\psi:\{0,1\}^{[2m]}\goesto[0,\infty)$ given by
$\psi(S)=\binom{m}{|S\cap[m]|}^{-1}\varphi(S\cap[m])$.  By
claim (2) above, $\psi^{\a}$ is NA for all $\a\in\RR^{2m}$ 
with $\a>\zero$.  By considering those $\a\in\RR^{2m}$ for which
$a_{i}=1$ for all $m+1\leq i\leq 2m$, it follows that $\varphi^{\a}$ is NA
for all $\a\in\RR^{m}$ with $\a>0$.
\end{proof}

\begin{coro}[Theorem 5.2 of \cite{BBL}]
If the initial distribution $\varphi_0$ of a SEP is such that $Z(\varphi)$
is stable (but not necessarily homogeneous), then
$Z(\varphi_t)$ is stable, and hence $\varphi_t$ is NA, for all $t\geq 0$.
\end{coro}

It is natural to try extending these results to asymmetric exclusion processes.
For $(i,j)\in \Lambda^{2}$ define $\pp_{ij}\in\EE$ by
$\pp_{ij}(S)=S\circ\tau_{ij}$ if $S(i)=1$ and $S(j)=0$,
and $\pp_{ij}(S)=S$ otherwise, for all $S\in\Omega$.
That is, $\pp_{ij}$ makes a particle hop from site $i$ to site $j$, if possible.
Let $E$ be a set of ordered pairs in $\Lambda^2$,
and for $(i,j)\in E$ let $\lambda_{ij}>0$.  An asymmetric exclusion process
is a Markov chain on $\Omega$ with semigroup $T(t)=\exp(t\L)$
generated by something of the form
$$
\L = \sum_{(i,j)\in E} \lambda_{ij} (\pp_{ij}-1).
$$

By the argument for claim (3) above, in order to show that $T(t)$
preserves stability for all $t\geq 0$, it suffices to do so for the two-site
semigroup $T_{\{1,2\}}(t)=\exp(t\L_{\{1,2\}})$ generated by
$$
\L_{\{1,2\}} = \lambda_{12}(\pp_{12}-1)+\lambda_{21}(\pp_{21}-1).
$$

\begin{exercise}[Strengthening Remark 5.3 of \cite{BBL}]
With the notation above, let $\lambda=\lambda_{12}+\lambda_{21}$,
$\beta_{12}=\lambda_{12}/\lambda$, and 
$\beta_{21}=\lambda_{21}/\lambda$.\\
\textup{(a)}\  In $\Alg$, $\pp_{12}+\pp_{21}=1+\tau_{12}$.\\
\textup{(b)}\  If $\omega$ is any word in $\{\pp_{12},\pp_{21}\}^n$,
then $\pp_{12}\omega=\pp_{12}$ and $\pp_{21}\omega=\pp_{21}$.\\
\textup{(c)}\  The semigroup generated by $\L_{\{1,2\}}$ is
$$
T_{\{1,2\}}(t) = \e^{-\lambda t} + (1-\e^{-\lambda t})
(\beta_{12}\pp_{12} + \beta_{21}\pp_{21})
$$
\textup{(d)}\  The semigroup $T_{\{1,2\}}(t)$ preserves stability for
all $t\geq 0$ if and only if $\beta_{12}=\beta_{21}=1/2$, in
which case it reduces to the SEP (of rate $\lambda/2$).\\
Thus, even the \underline{slightest} asymmetry ruins preservation
of stability by the SEP!
\end{exercise}

Finally, we consider a SEP in which particle number is not conserved.
For $i\in\Lambda$ define $\aa_{i},\aa_{i}^{*}\in
\EE$ as follows:\ for $S\in\Omega$ and 
$j\in\Lambda$, let $(\aa_{i}(S))(j)=(\aa_{i}^{*}(S))(j)=S(j)$ if 
$j\neq i$, and $(\aa_{i}(S))(i)=0$ and $(\aa_{i}^{*}(S))(i)=1$.
That is, $\aa_{i}$ annihilates a particle at site $i$, and
$\aa_{i}^{*}$ creates a particle at site $i$, if possible.

A SEP with particle creation and annihilation is a Markov chain on
$\Omega$ with semigroup $T(t)=\exp(t\L)$ generated by something of
the form
$$
\L = \sum_{\{i,j\}\in E} \lambda_{ij} (\tau_{ij}-1) +
\sum_{i\in \Lambda}
\left[ \theta_{i}(\aa_{i}-1) + \theta_{i}^{*}(\aa_{i}^{*}-1) \right],
$$
in which the first sum is the generator of the SEP in Theorem 7.1 and
$\theta_{i},\theta_{i}^{*}\geq 0$ for each $i\in\Lambda$.

By the argument for claim (3) above, to show that this $T(t)$
preserves stability for all $t\geq 0$, it suffices to do so for the
one-site semigroups generated by $\L_{1} = \theta(\aa_{1}-1)$ and
$\L_{1}^{*} = \theta(\aa_{1}^{*}-1)$, respectively.

\begin{exercise}
The semigroups generated by $\L_1$ and $\L_1^*$ are
$$
T_{1}(t) = \e^{-\theta t} + (1-\e^{-\theta t})\aa_{1}
\ \ \ \mathrm{and}\ \ \ 
T_{1}^{*}(t) = \e^{-\theta t} + (1-\e^{-\theta t})\aa_{1}^{*},
$$
respectively.  Both $T_{1}(t)$ and $T_{1}^{*}(t)$ preserve stability.
\end{exercise}

\begin{coro}
If the initial distribution $\varphi_0$ of a SEP with particle
creation and annihilation is such that $Z(\varphi)$ is stable, then
$Z(\varphi_t)$ is stable, and hence $\varphi_t$ is NA, for all $t\geq 0$.
\end{coro}

\section{Inequalities for mixed discriminants.}

This section summarizes a powerful application of stable polynomials from
a 2008 paper of Gurvits \cite{Gu}.

We will use without mention the facts that $\log$ and $\exp$ are strictly
increasing functions on $(0,\infty)$.
A function $\rho:I\goesto\RR$ defined on an interval $I\subseteq\RR$ is
\emph{convex} provided that for all $a_1,a_2\in I$,
$\rho((a_1+a_2)/2)\leq (\rho(a_1)+\rho(a_2))/2$.  It is \emph{strictly
convex} if it is convex and equality holds here only when $a_1=a_2$.
A function $\rho:I\goesto\RR$ is \emph{(strictly) concave} if $-\rho$
is (strictly) convex.
For example, for positive reals $a_1,a_2>0$ one has
$(\sqrt{a_1}-\sqrt{a_2})^2\geq0$, with equality only if $a_1=a_2$.
It follows that $\log((a_1+a_2)/2)\geq (\log(a_1) + \log(a_2))/2$, with
equality only if $a_1=a_2$.  That is, $\log$ is strictly concave.\\

\noindent\textbf{Jensen's Inequality} (Theorem 90 of \cite{HLP})\textbf{.}\
\emph{
Let $\rho:I\goesto\RR$ be defined on an interval $I\subseteq\RR$,
let $a_i\in I$ for $i\in[n]$, and let $b_i>0$ for $i\in[n]$ be such
that $\sum_{i=1}^n b_i = 1$.  If $\rho$ is convex then
$$
\rho\left(\sum_{i=1}^n b_i a_i\right) \leq
\sum_{i=1}^n b_i\, \rho(a_i).
$$
If $\rho$ is strictly convex and equality holds, then $a_1=a_2=\cdots=a_n$.}\\

For integer $d\geq 1$, let $\G(d) = (1-1/d)^{d-1}$, and let $\G(0)=1$.
Note that $\G(1)=0^0=1$, and that $\G(d)$ is a strictly decreasing function
for $d\geq 1$.  For homogeneous $f\in\RR[\x]$ with nonnegative coefficients,
define the \emph{capacity} of $f$ to be
$$
\CAP(f) = \inf_{\c>\zero} \frac{f(\c)}{c_1\cdots c_m},
$$
with the infimum over the set of all $\c\in\RR^m$ with $c_i>0$ for all $i\in[m]$.

\begin{lemma}[Lemma 3.2 of \cite{Gu}]
Let $f=\sum_{i=0}^d b_i x\in\RR[x]$ be a nonzero
univariate polynomial of degree $d$ with nonnegative coefficients.
If $f$ is real stable then
$b_1=f'(0)\geq \G(d)\CAP(f)$, and if $\CAP(f)>0$ then
equality holds if and only if $d\leq 1$ or $f(x) = b_d(x+\xi)^d$ for some $\xi>0$.
\end{lemma}
\begin{proof}
If $\CAP(f)=0$ then there is nothing to prove,
so assume that $\CAP(f)>0$.
If $d=0$ then $f'(0)= b_{1} = 0 = \G(0)\CAP(f)$, and
if $d=1$ then $f'(0)= b_1 = \G(1)\CAP(f)$, so assume that $d\geq 2$.
If $f(0)=0$ then $f'(0)=\lim_{c\goesto 0}f(c)/c\geq \CAP(f) > 
\G(d)\CAP(f)$.
Thus, assume that $d\geq 2$ and $f(0)=b_0>0$.  We may rescale the polynomial
so that $b_0=1$.  Now there are $a_i>0$ for $i\in[d]$ such that
$$
f(x) = \prod_{i=1}^d (1+a_i x),
$$
and $b_1=a_1+\cdots+a_d$.  For any $c>0$ we have
$$
\frac{\log(\CAP(f)c)}{d}
\leq
\frac{\log(f(c))}{d} = \frac{1}{d}\sum_{i=1}^d \log(1+a_i c)
\leq
\log\left(1+\frac{b_1 c}{d}\right),
$$
by Jensen's Inequality.  It follows that $\CAP(f)c\leq (1+b_1c/d)^d$
for all $c>0$.  Let $g(x)=(1+b_1x/d)^d$.  Elementary calculus shows that
$$
\CAP(g)= \inf_{c>0} \frac{g(c)}{c} = \frac{g(c_*)}{c_*} = \frac{b_1}{\G(d)},
\ \ \mathrm{in\ which}\ \ c_*=\frac{d}{b_1(d-1)}.
$$
Since $\CAP(f)\leq \CAP(g)$, this yields the stated inequality.
If equality holds, then equality holds in the application
of Jensen's Inequality, and so $f$ has the stated form.
\end{proof}

\begin{lemma}[Theorem 4.10 of \cite{Gu}]
Let $f\in\StabR[x_1,...,x_m]$ be real stable, with nonnegative 
coefficients, and homogeneous of degree $m$.
Let $g = \partial_m f|_{x_m=0}$.  Then
$$
\CAP(g) \geq \G(\deg_m(d))\CAP(f).
$$
\end{lemma}
\begin{proof}
We may assume that $d=\deg_m(f)\geq 1$.
Let $c_i>0$ for $i\in[m-1]$, and let $p_\c(x)=f(c_1,...,c_{m-1},x)$.
Since $f$ has nonnegative coefficients, $p_{\c}\not\equiv 0$.
As in the proof of Lemma 2.4(f), $p_\c$ has degree $d$.
By specialization, $p_\c$ is real stable.  Lemma 8.1 implies that
$$
g(\c)=p_\c'(0)\geq \G(d)\CAP(p_\c)\geq\G(d)\CAP(f)
$$
for all $\c\in\RR^{m-1}$ with $\c>\zero$.  If $m=1$ then $g=\CAP(g)$
is a constant.  If $m\geq 2$ then for any such $\c$ let
$b=(c_1\cdots c_{m-1})^{-1/(m-1)}$.  Since $g$ is homogeneous of
degree $m-1$,
$$
\frac{g(\c)}{c_1\cdots c_{m-1}}=g(bc_1,...,bc_{m-1})
\geq 
\G(d)\CAP(f).
$$
It follows that $\CAP(g)\geq \G(d)\CAP(f)$.
\end{proof}

\begin{theorem}[Theorem 2.4 of \cite{Gu}]
Let $f\in\StabR[x_1,...,x_m]$ be real stable, with nonnegative 
coefficients, and homogeneous of degree $m$.  Let $\deg_{i}(f)=d_i$
and $e_i=\min\{i,d_i\}$ for each $i\in[m]$.  Then
$$
\del^{\one} f(\zero) \geq \CAP(f)\prod_{i=2}^m \G(e_i).
$$
\end{theorem}
\begin{proof}
Let $g_m=f$ and let $g_{i-1}=\partial_i g_i|_{x_i=0}$ for all
$i\in[m]$.  By contraction and specialization, $g_{i}$ is real stable
for each $i\in[m]$.  Notice that $g_0=\del^{\one} f(\zero)=\CAP(g_0)$.
By Lemma 8.2.
$\CAP(g_{i-1})\geq\CAP(g_i)\cdot \G(\deg_i g_i)$ for each $i\in[m]$.
But $\deg_i g_i \leq \deg_i f = d_i$, and $\deg_i g_i$ is at most the
total degree of $g_i$, which is $i$.  Hence $\deg_i g_i \leq e_i$,
and thus $\G(\deg_i g_i) \geq \G(e_i)$.  Thus
$\CAP(g_{i-1})\geq\CAP(g_i)\cdot \G(e_{i})$ for each $i\in[m]$.
Combining these inequalities (and $\G(e_1)=1$) gives the result.
\end{proof}

With the notation of Theorem 8.3, since $e_i\leq i$ for all $i\in[m]$
and $G(d)$ is a decreasing function of $d$, one has the inequality
$$
\prod_{i=2}^m \G(e_i) \geq \prod_{i=2}^m \G(i)
= \prod_{i=2}^m \left(\frac{i-1}{i}\right)^{i-1} = \frac{m!}{m^m}.
$$
Thus, the following corollary is immediate.
\begin{coro}
Let $f\in\StabR[x_1,...,x_m]$ be real stable, with nonnegative 
coefficients, and homogeneous of degree $m$.  Then
$$
\del^{\one} f(\zero) \geq \frac{m!}{m^m}\cdot\CAP(f).
$$
\end{coro}

\begin{theorem}[Theorem 5.7 of \cite{Gu}]
Let $f\in\StabR[x_1,...,x_m]$ be real stable, with nonnegative 
coefficients, and homogeneous of degree $m$.  Equality holds in the bound of
Corollary 8.4 if and only if there are nonnegative reals $a_i\geq 0$ for
$i\in[m]$ such that
$$
f(\x) = (a_1 x_1 + \cdots + a_m x_m)^m.
$$
\end{theorem}
(We omit the proof.)

\begin{lemma}[Fact 2.2 of \cite{Gu}]
Let $f\in\RR[x_1,...,x_m]$ be homogeneous of degree $m$, with nonnegative
coefficients.  Assume that $\partial_{i}f(\one)=1$ for all $i\in[m]$.  Then $\CAP(f)=1$.
\end{lemma}
\begin{proof}
Let $f=\sum_{\alpha} b(\alpha) \x^{\alpha}$, so that if 
$b(\alpha)\neq 0$ then $|\alpha|=\sum_{i=1}^{m}\alpha(i)=m$.
By hypothesis, for all $i\in[m]$, $\sum_{\alpha} b(\alpha)\alpha(i)=1$.
Averaging these over all $i\in[m]$ yields $f(\one)=\sum_{\alpha} b(\alpha)=1$,
so that $\CAP(f)\leq 1$.  Conversely, let $\c\in\RR^{m}$ with $\c>\zero$.
Jensen's Inequality implies that
\showon
\log(f(\c))
&=&
\log\left(\sum_{\alpha} b(\alpha) \c^{\alpha}\right)\\
&\geq&
\sum_{\alpha} b(\alpha) \log(\c^{\alpha})
=
\sum_{i=1}^{m}\log(c_{i})\sum_{\alpha} b(\alpha)\alpha(i)
=
\log( c_{1}\cdots c_{m}).
\showoff
It follows that $\CAP(f)\geq 1$.
\end{proof}

\begin{example}[van der Waerden Conjecture]
An $m$-by-$m$ matrix $A=(a_{ij})$ is \emph{doubly stochastic} if 
all entries are nonnegative reals and every row and column sums to one.
In 1926, van der Waerden conjectured that if $A$ is an $m$-by-$m$ doubly
stochastic matrix then $\per(A) \geq m!/m^m$, with equality if and only if
$A=(1/m)J$, the $m$-by-$m$ matrix in which every entry is $1/m$.
In 1981 this lower bound was proved by Falikman, and the characterization of
equality was proved by Egorychev.  These results follow immediately from
Corollary 8.4 and Theorem 8.5, as follows.  It suffices to prove the
result for an $m$-by-$m$ doubly stochastic matrix $A=(a_{ij})$ with
no zero entries, by a routine limit argument.  The polynomial
$$
f_A(\x) = \prod_{j=1}^m\left(a_{1j} x_1 + \cdots + a_{mj} x_m\right)
$$
is clearly homogeneous and real stable, with nonnegative coefficients and
of degree $m$, and such that $\deg_{i}(f_{A})=m$ for all $i\in[m]$.   Since $A$
is doubly stochastic, Lemma 8.6 implies that $\CAP(f_{A})=1$.
Since
$$
\per(A) = \del^{\one} f_A(\zero),
$$
Corollary 8.4 and Theorem 8.5 imply the results of Falikman and Egorychev,
respectively.  Gurvits \cite{Gu} also uses a similar argument to prove a refinement
of the van der Waerden conjecture due to Schrijver and Valiant -- see 
also \cite{LS}.
\end{example}

Given $n$-by-$n$ matrices $A_1$,...,$A_m$, the \emph{mixed discriminant} of
$\AA = (A_1,...,A_m)$ is
$$
\Disc(\AA) =
\left. \del^{\one} \det(x_1 A_1 + \cdots + x_m A_m)\right|_{\x=\zero}.
$$
This generalizes the permanent of an $m$-by-$m$ matrix $B=(b_{ij})$
by considering the collection of matrices $\AA(B) = (A_1,...,A_m)$
defined by $A_h = \diag(a_{h1},...,a_{hm})$ for each $h\in[m]$.
In this case one sees that 
$$
\det(x_1 A_1 + \cdots + x_m A_m) = f_B(\x)
$$
with the notation of Example 8.7, and it follows that
$\Disc(\AA(B)) = \per(B)$.

\begin{example}[Bapat's Conjecture]
Generalizing the van der Waerden conjecture, in 1989 Bapat considered
the set $\Omega(m)$ of $m$-tuples of $m$-by-$m$ matrices
$\AA = (A_1,...,A_m)$ such that each $A_i$ is positive semidefinite
with trace $\tr(A_i)=1$, and $\sum_{i=1}^m A_i = I$.  For any doubly stochastic
matrix $B$, $\AA(B)$ is in this set.  The natural conjecture is that for
all $\AA\in\Omega(m)$, $\Disc(\AA)\geq m!/m^m$, and equality is
attained if and only if $\AA=\AA((1/m)J)$.  This was proved by Gurvits in 2006 --
again, it follows directly from Corollary 8.4 and Theorem 8.5.
It suffices to prove the result for $\AA\in\Omega(m)$ such that
each $A_{i}$ is positive definite, by a routine limit argument. 
By Proposition 2.1, for $\AA\in\Omega(m)$, the polynomial
$$
f_\AA(\x) = \det(x_1 A_1 + \cdots + x_m A_m)
$$
is real stable.  Since each $A_i$ is positive definite, all coefficients
of $f_\AA$ are nonnegative, $f_\AA$ is homogeneous of degree 
$m$,  and $\deg_{i}(f_{A})=m$ for all $i\in[m]$.
Since $\AA\in \Omega(m)$, Lemma 8.6 implies that 
$\CAP(f_{\AA})=1$.  Thus, $f_\AA$ satisfies the hypothesis of Theorems 
8.3 and 8.5, and since $\Disc(\AA) = \del^{\one} f_\AA(\zero)$,
the result follows.
\end{example}

\section{Further Directions.}

\subsection{Other circular regions.}

Let $\Omega\subseteq\CC^m$.  A polynomial $f\in\CC[\x]$ is
\emph{$\Omega$-stable} if either $f\equiv 0$ identically, or
$f(\z)\neq 0$ for all $\z\in\Omega$.  At this level of generality little
can be said.  If $\Omega=\A_1\times\cdots\times\A_m$ is a product of
open circular regions then there are M\"obius transformations
$z\mapsto \phi_i(z)=(a_i z+b_i)/(c_i z+d_i)$ such that $\phi_i(\HH)=\A_i$
for all $i\in[m]$.  The argument in Section 4.1 shows that
$f\in\CC[\x]$ is $\Omega$-stable if and only if
$$
\widetilde{f}=(\c\z+\mathbf{d})^{\deg f}\cdot f(\phi_1(z_1),...,\phi_m(z_m))
$$
is stable.  In this way results about stable polynomials can be translated
into results about $\Omega$-stable polynomials for any $\Omega$ that is
a product of open circular regions.

Theorem 6.3 of \cite{BB5} is the $\Omega$-stability analogue of Theorem 5.2.
We mention only two consequences of this.  Let $\D=\{z\in\CC:\ |z|<1\}$ be the
open unit disc, and for $\theta\in\RR$ let $\HH_\theta=\{\e^{-\i\theta}z:\
z\in\HH\}$.  Thus $\HH_0=\HH$, and $\HH_{\pi/2}$ is the open right half-plane.
A $\D^m$-stable polynomial is called \emph{Schur stable}, and a
$\HH_{\pi/2}^m$-stable polynomial is called \emph{Hurwitz stable}.

\begin{prop}[Remark 6.1 of \cite{BB5}.]
Fix $\kappa\in\NN^m$, and let $T:\CC[\x]^{\leq\kappa}\goesto\CC[\x]$
be a linear transformation.  The following are equivalent:\\
\textup{(a)}\ $T$ preserves Schur stability.\\
\textup{(b)}\ $T((\one+\x\y)^\kappa)$ is Schur stable in $\CC[\x,\y]$.
\end{prop}

\begin{prop}[Remark 6.1 of \cite{BB5}.]
Fix $\kappa\in\NN^m$, and let $T:\CC[\x]^{\leq\kappa}\goesto\CC[\x]$
be a linear transformation.  The following are equivalent:\\
\textup{(a)}\ $T$ preserves Hurwitz stability.\\
\textup{(b)}\ $T((\one+\x\y)^\kappa)$ is Hurwitz stable in $\CC[\x,\y]$.
\end{prop}

\subsection{Applications of Theorem 5.4.}

It is natural to consider a multivariate anaogue of the multiplier sequences
studied by P\'olya and Schur.  Let $\lambda:\NN^m\goesto\RR$, and define
a linear transformation $T_\lambda:\CC[\x]\goesto\CC[\x]$ by
$T_\lambda(\x^\alpha)= \lambda(\alpha)\x^\alpha$ for all $\alpha\in\NN^m$,
and linear extension.  For which $\lambda$ does $T_{\lambda}$
preserve real stability?  The answer:\  just the ones you get from
the P\'olya-Schur Theorem, and no more.

\begin{theorem}[Theorem 1.8 of \cite{BB4}.]
Let $\lambda:\NN^m\goesto\RR$.  Then $T_\lambda$ preserves real stability if
and only if there are univariate multiplier sequences $\lambda_i:\NN\goesto\RR$
for $i\in[m]$ and $\epsilon\in\{-1,+1\}$ such that
$$
\lambda(\alpha) = \lambda_1(\alpha(1))\cdots\lambda_m(\alpha(m))
$$
for all $\alpha\in\NN^m$, and either
$\epsilon^{|\alpha|}\lambda(\alpha)\geq 0$ for all $\alpha\in\NN^m$, or
$\epsilon^{|\alpha|}\lambda(\alpha)\leq 0$ for all $\alpha\in\NN^m$.
\end{theorem}

Theorem 5.4 (and similarly Propositions 9.1 and 9.2) can be used to derive
a wide variety of results of the form:\ such-and-such an operation
preserves stability (or Schur or Hurwitz stability).  Here is a short
account of Hinkkanen's proof of the Lee-Yang Circle Theorem, taken
from Section 8 of \cite{BB6}.

For $f,g\in\CC[\x]^\ma$, say $f=\sum_{S\subseteq[m]} a(S)\x^S$ and
$g=\sum_{S\subseteq[m]} b(S)\x^S$, let
$$
f\bullet g = \sum_{S\subseteq[m]} a(S)b(S) \x^S
$$
be the \emph{Schur-Hadamard product} of $f$ and $g$.

\begin{theorem}[Hinkkanen, Theorem 8.5 of \cite{BB6}]
If $f,g\in\CC[\x]^\ma$ are Schur stable then $f\bullet g$ is Schur
stable.
\end{theorem}
\begin{proof}
Let $T_g:\CC[\x]^\ma\goesto\CC[\x]^\ma$ be defined by $f\mapsto f\bullet g$.
By Proposition 9.1, to show that
$T_g$ preserves Schur stability it suffices to show that $T_g((\one+\x\y)^{[m]})$
is Schur stable.  Clearly $T_g((\one+\x\y)^{[m]}) = g(x_1 y_1,...,x_m y_m)$
is Schur stable since $g(\x)$ is.  Hence $T_g$ preserves Schur stability,
and so $f\bullet g$ is Schur stable.
\end{proof}

\begin{theorem}[Lee-Yang Circle Theorem, Theorem 8.4 of \cite{BB6}]
Let $A=(a_{ij})$ be a Hermitian $m$-by-$m$ matrix with $|a_{ij}|\leq 1$
for all $i,j\in[m]$.  Then the polynomial
$$
f(\x) = \sum_{S\subseteq[m]}\x^S \prod_{i\in S}\prod_{j\not\in S} a_{ij}
$$
is Schur stable.  The diagonalization $g(x)=f(x,...,x)$ is such that
$x^m g(1/x) = g(x)$, and it follows that all roots of
$g(x)$ are on the unit circle.
\end{theorem}
\begin{proof}
For $i<j$ in $[m]$ let
$$
f_{ij} = (1+a_{ij} x_i + \overline{a_{ij}} x_j + x_i x_j)
\prod_{h\in[m]\drop\{i,j\}} (1+x_h).
$$
One can check that each $f_{ij}$ is Schur stable.  The polynomial $f(\x)$
is the Schur-Hadamard product of all the
$f_{ij}$ for $\{i,j\}\subseteq[m]$.  By Theorem 9.4, $f(\x)$ is Schur stable.
\end{proof}

Section 8 of \cite{BB6} contains many many more results of this nature.

\subsection{A converse to the Grace-Walsh-Szeg\H{o} Theorem.}

The argument of Sections 4.2 and 4.3 can be used to prove the following.

\begin{exercise}
If $f\in\Stab[\x]^\ma$ is multiaffine and stable then
$$
T_{\S(m)} (f) = \frac{1}{m!} \sum_{\sigma\in\S(m)} \sigma(f)
$$
is multiaffine and stable.
\end{exercise}

This is in fact equivalent to the GWS Theorem,
since for all $f\in\CC[\x]^{\ma}$, $T_{\S(m)}f(\x) = \Pol_{m} f(x,\ldots,x)$.
For which transitive permutation groups $G\leq\S(m)$ does
the linear transformation $T_G = |G|^{-1}\sum_{\sigma\in G}\sigma$
preserve stability?  The answer:\ not many, and they give nothing new.

\begin{theorem}[Theorem 6 of \cite{BW}.]
Let $G\leq\S(m)$ be a transitive permutation group such that
$T_G$ preserves stability.  Then $T_G=T_{\S(m)}$.
\end{theorem}

\subsection{Phase and support theorems.}

A polynomial $f\in\CC[\x]$ has \emph{definite parity} if every monomial
$\x^\alpha$ occurring in $f$ has total degree of the same parity:\
all are even, or all are odd.

\begin{theorem}[Theorem 6.2 of \cite{COSW}]
Let $f\in\CC[\x]$ be Hurwitz stable and with definite parity.  Then
there is a phase $0\leq \theta < 2\pi$ such that $\mathrm{e}^{-\i\theta}f(\x)$
has only real nonnegative coefficients.
\end{theorem}

The \emph{support} of $f=\sum_{\alpha}c(\alpha)\x^{\alpha}$ is
$\supp(f) = \{\alpha\in\NN^{m}:\ c(\alpha)\neq 0\}$.
Let $\delta_i$ denote the unit vector with a one in the $i$-th coordinate,
and for $\alpha\in\ZZ^{n}$ let $|\alpha|=\sum_{i=1}^m |\alpha(i)|$.
A \emph{jump system} is a subset $\J\subseteq\ZZ^m$ satisfying the following
\emph{two-step axiom}:\\
\textbf{(J)}\  If $\alpha,\beta\in\J$ and $i\in[m]$ and
$\epsilon\in\{-1,+1\}$ are such that $\alpha'=\alpha+\epsilon\delta_i$
satisfies $|\alpha'-\beta|<|\alpha-\beta|$, then either $\alpha'\in\J$
or there exists $j\in[m]$ and $\varepsilon\in\{-1,+1\}$ such that
$\alpha''=\alpha'+\varepsilon\delta_j\in\J$ and
$|\alpha''-\beta|<|\alpha'-\beta|$.

Jump systems generalize some more familiar combinatorial objects.
A jump system contained in $\{0,1\}^{m}$ is a \emph{delta-matroid}.
A delta-matroid $\J$ for which $|\alpha|$ is constant for all $\alpha\in\J$
is the set of bases of a \emph{matroid}.  For bases of matroids, the
two-step axiom (J) reduces to the basis exchange axiom familiar from linear 
algebra:\ if $A,B\in\J$ and $a\in A\drop B$, then there exists $b\in 
B\drop A$ such that $(A\drop\{a\})\cup\{b\}$ is in $\J$.

\begin{theorem}[Theorem 3.2 of \cite{B}]
If $f\in\Stab[\x]$ is stable then the support $\supp(f)$ is a jump system.
\end{theorem}

Recall from Section 7 that for multiaffine polynomials with nonnegative
coefficients, real stability implies the Rayleigh property.
A set system $\J$ is \emph{convex} when $A,B\in\J$ and $A\subseteq B$ imply 
that $C\in\J$ for all $A\subseteq C\subseteq B$.
\begin{theorem}[Section 4 of \cite{W}]
Let $f=\sum_{S\subseteq[m]} c(S) \x^S$ be multiaffine with real nonnegative
coefficients, and assume that $f$ is Rayleigh.\\
\textup{(a)}\ The support $\supp(f)$ is a convex delta-matroid.\\
\textup{(b)}\ The coefficients are \emph{log-submodular}:\ for all
$A,B\subseteq[m]$,
$$
c(A\cap B) c(A\cup B) \leq c(A) c(B).
$$
\end{theorem}

\bibliographystyle{amsplain}

\end{document}